\documentclass[onefignum,onetabnum]{siamart220329}
\graphicspath{ {Figures/} }
\usepackage{makecell}
\usepackage{caption}
\usepackage{subcaption,microtype}

\usepackage{lipsum}
\usepackage{amsfonts}
\usepackage{graphicx}
\usepackage{epstopdf}
\usepackage{algorithmic}
\ifpdf
  \DeclareGraphicsExtensions{.eps,.pdf,.png,.jpg}
\else
  \DeclareGraphicsExtensions{.eps}
\fi

\newsiamremark{remark}{Remark}
\newsiamremark{hypothesis}{Hypothesis}
\crefname{hypothesis}{Hypothesis}{Hypotheses}
\newsiamthm{claim}{Claim}

\headers{Mathematical modelling of pesticide leaf uptake}{J. Delos Reyes, T. Shardlow, M.B. Delgado-Charro, S. Webb, K.A.J. White}

\title{Mathematical modelling of adjuvant-enhanced active ingredient leaf uptake of pesticides\thanks{
\textbf{Funding: }{This work was funded by Engineering and Physical Sciences Research Council grant EP/S515279/1 and by Syngenta UK Ltd. grant TK0448301.}}}

\author{ Jenny Delos Reyes\thanks{Corresponding author. Department of Mathematical Sciences, University of Bath, Bath, UK, BA2 7AY 
    (\email{jdr47@bath.ac.uk}).}
\and Tony Shardlow\thanks{Department of Mathematical Sciences, University of Bath, Bath, UK, BA2 7AY.}
\and M. Bego\~{n}a Delgado-Charro\thanks{Department of Life Sciences, University of Bath, Bath, UK, BA2 7AY.}
\and Steven Webb\thanks{Product Safety, Syngenta, Jealott's Hill International Research Centre, Bracknell, UK, RG42 6EY.}
\and K. A. Jane White\footnotemark[3]}

\usepackage{amsopn}

\ifpdf
\hypersetup{
  pdftitle={Mathematical modelling of adjuvant-enhanced active ingredient leaf uptake of pesticides},
  pdfauthor={J. Delos Reyes, T. Shardlow, M. B. Delgado-Charro, S. Webb, K. A. J. White}
}
\fi

\begin{document}

    \maketitle

    % REQUIRED
    \begin{abstract}
        The global importance of effective and affordable pesticides to optimise crop yield and to support health of our growing population cannot be understated. But to develop new products or refine existing ones in response to climate and environmental changes is both time-intensive and expensive which is why the agrochemical industry is increasingly interested in using mechanistic models as part of their formulation development toolbox. In this work, we develop such a model to describe uptake of pesticide spray droplets across the leaf surface. We simplify the leaf structure by identifying the outer cuticle as the main barrier to uptake; the result is a novel, hybrid model in which two well-mixed compartments are separated by a membrane in which we describe the spatio-temporal distribution of the pesticide. This leads to a boundary value partial differential equation problem coupled to a pair of ordinary differential equation systems which we solve numerically. We also simplify the pesticide formulation into two key components: the Active Ingredient which produces the desired effect of the pesticide and an Adjuvant which is present in the formulation to facilitate effective absorption of the Active Ingredient into the leaf. This approach gives rise to concentration-dependent diffusion. We take an intuitive approach to parameter estimation using a small experimental data set and subsequently demonstrate the importance of the concentration-dependent diffusion in replicating the data. Finally, we demonstrate the need for further work to identify how the physicochemical properties of pesticides affect flow into and across the leaf surface. 
    \end{abstract}

    % REQUIRED

    % REQUIRED

    \section{Introduction}
	\label{sec: Introduction}
        The importance of effective and affordable agrochemicals worldwide cannot be understated \cite{YeungMayT2017DICo}. They are used globally to optimise crop yield in a number of different ways including growth enhancement (for example, herbicides that kill unwanted plants to eliminate competition for resources \cite{FedtkeCarl1982Bapo,ArandKatja2018TMoA}) and disease management (for example, organophosphate insecticides which kill mosquitoes to control spread of diseases such as malaria and dengue \cite{RiveroAna2010Icov}). Traditionally the process of developing a new product is lengthy and expensive, and this has helped to strengthen recent interest in adding mathematical models to the product development toolbox \cite{KalyabinaValeriyaP2021Pfdp}. The model which we present here is our contribution and we focus on the uptake of pesticides through the leaf surface. 

        There is a small literature on mechanistic models for pesticide uptake in plants as summarised in \cite{FANTKE20111639}, but often these models focus on root exposure (see, for example, \cite{TrappStefan1995GOMf,Fujisawa2002}) and not direct leaf surface contact with the pesticide spray. An uptake model that incorporates foliar exposure and whole plant allocation of absorbed chemical has been presented in \cite{SatchiviNorbertM2000ANDSa,SatchiviNorbertM2000ANDSb}, and was used here as a basis for model development.
        
\begin{figure}[H]
		\includegraphics[width=10cm]{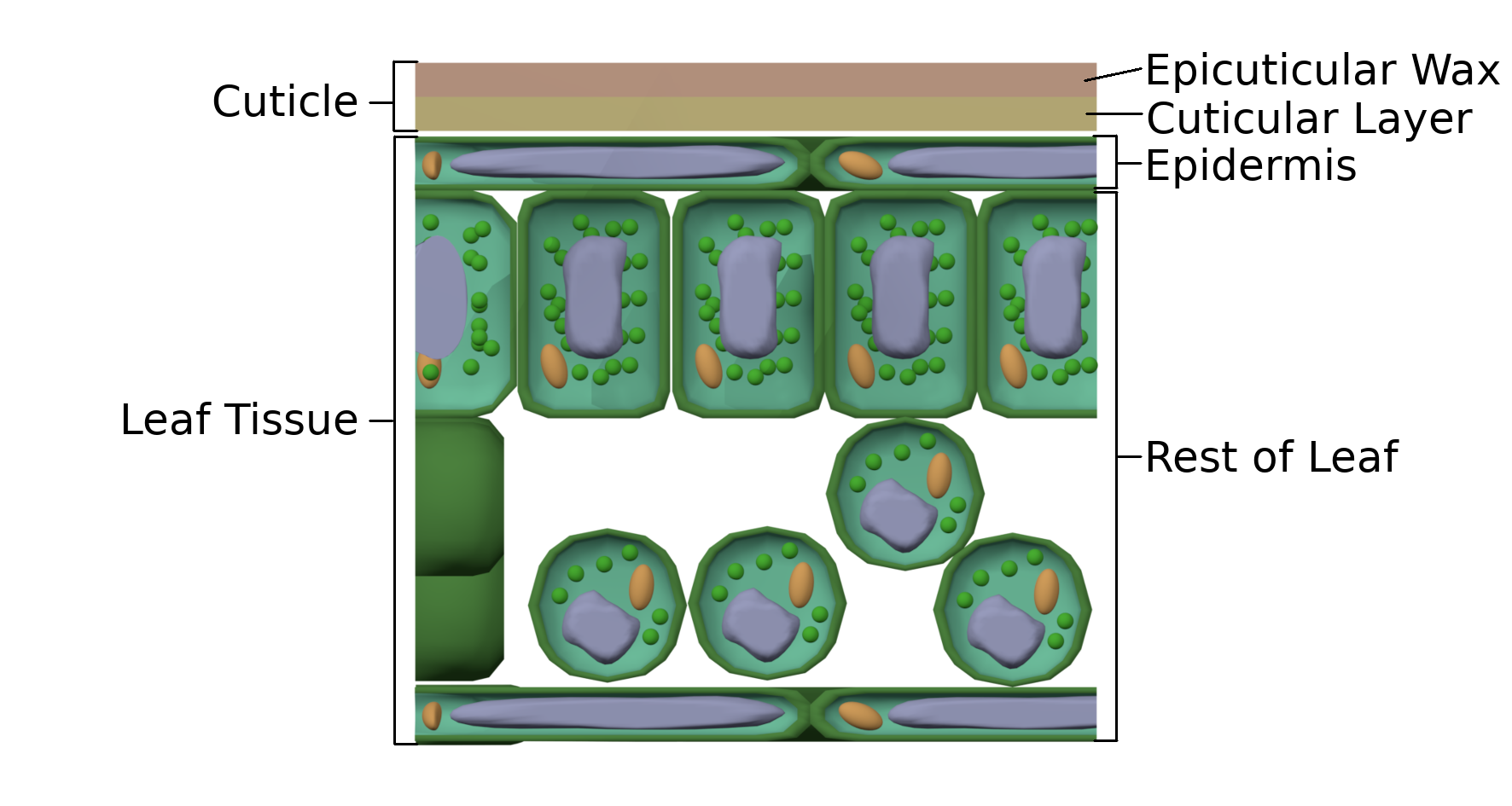}
		\centering
		\caption{Simplified diagram of the leaf structure.}
		\label{fig: leaf_model}
	\end{figure}

        One challenge in modelling uptake across the leaf surface is to determine the appropriate simplifying assumptions about the leaf structure (see \cref{fig: leaf_model}) which consists of many layers and which is also laterally heterogeneous \cite{SatchiviNorbertM2000ANDSa,SchreiberLukas2009WaSP,LiShuai2013Thas}. One layer that is clearly identified as important in the uptake process is the uppermost layer, known as the cuticle. This layer is non-cellular, made of cutin which is a waxy substance \cite{Jeffree2006,SchreiberLukas2009WaSP}. The cuticle acts as the main barrier to prevent water loss or cuticular transpiration, and it protects the leaf from any external threat such as chemical attacks and ultraviolet radiation \cite{MulroyThomasW1979SPoH}.

        As with the leaf, the composition of pesticides are highly complex. Each formulation will necessarily include an {\bf Active Ingredient (AI)} because this is the compound that actually works against pests by controlling, killing, or repelling them. Beyond that, formulations include a wide range of other compounds each of which performs some function. One of the more common compounds in a formulation is known more generally as an {\bf Adjuvant (AJ)} which is added to facilitate the absorption of AI into the leaf, however the exact mechanisms are extremely complex, and only partially understood \cite{WangC.J2007Fuop, ArandKatja2018TMoA}. For example, an accelerator AJ \cite{Schonherr1993} will increase the rate at which the AI can move across the leaf cuticle \cite{SchreiberLukas2009WaSP} by increasing the flexibility or fluidity of the waxes and cutin \cite{SchonherrBaur1994}. Another common AJ is a surfactant \cite{surfactants-in-agrochemicals} which acts on the contact between pesticide droplet and leaf surface in order to optimise delivery of the AI. 

        With leaf structure and pesticide composition in mind, %whilst also exploiting Occam's Razor \cite{VanDenBergHugoA.2018OrfO}, 
        in the following \Cref{sec: Model} we formulate our model and elucidate our simple approach to model parameter estimation using a small data set. \Cref{sec: Results} presents our results which consist of numerical solutions of the model system. These numerics %are not trivial since they 
        involve solving a boundary value partial differential equation (PDE). Our results highlight the importance of the AJ in creating predictions that are consistent with experimental data. Finally in \Cref{sec: Discussion} we discuss limitations of our model and identify future work to develop more robust empirical relations between key model processes and the underlying physico-chemical properties of the pesticide components.

    \section{Model formulation}
	\label{sec: Model}
        In our model, we assume a simplified leaf anatomy comprising two layers, namely the wax/cuticle and the leaf tissue, together with a pesticide droplet which instantaneously settles to a stable configuration on the leaf surface. We consider the case of no evaporation, which means that all volumes and contact areas between layers and the pesticide droplet remain constant throughout. For simplicity we assume that the droplet shape can be represented by a hemisphere and that each droplet is identical and spatially segregated from all other droplets. We assume that the pesticide moves from the droplet across the leaf cuticle into the leaf tissue in the transverse direction only so that our problem becomes spatially one-dimensional. This approach, in which lateral movement within the cuticle is ignored, has been used widely in the literature for flow across biological membranes (see, for example, \cite{CleekRL1993Anmf,JeppsOwenG2013Mths,HansenSteffi2013Iipf,TodoHiroaki2013Mmtp,JonesJ.G2016Amat}).  

\begin{figure}[H]
		\includegraphics[width=13cm]{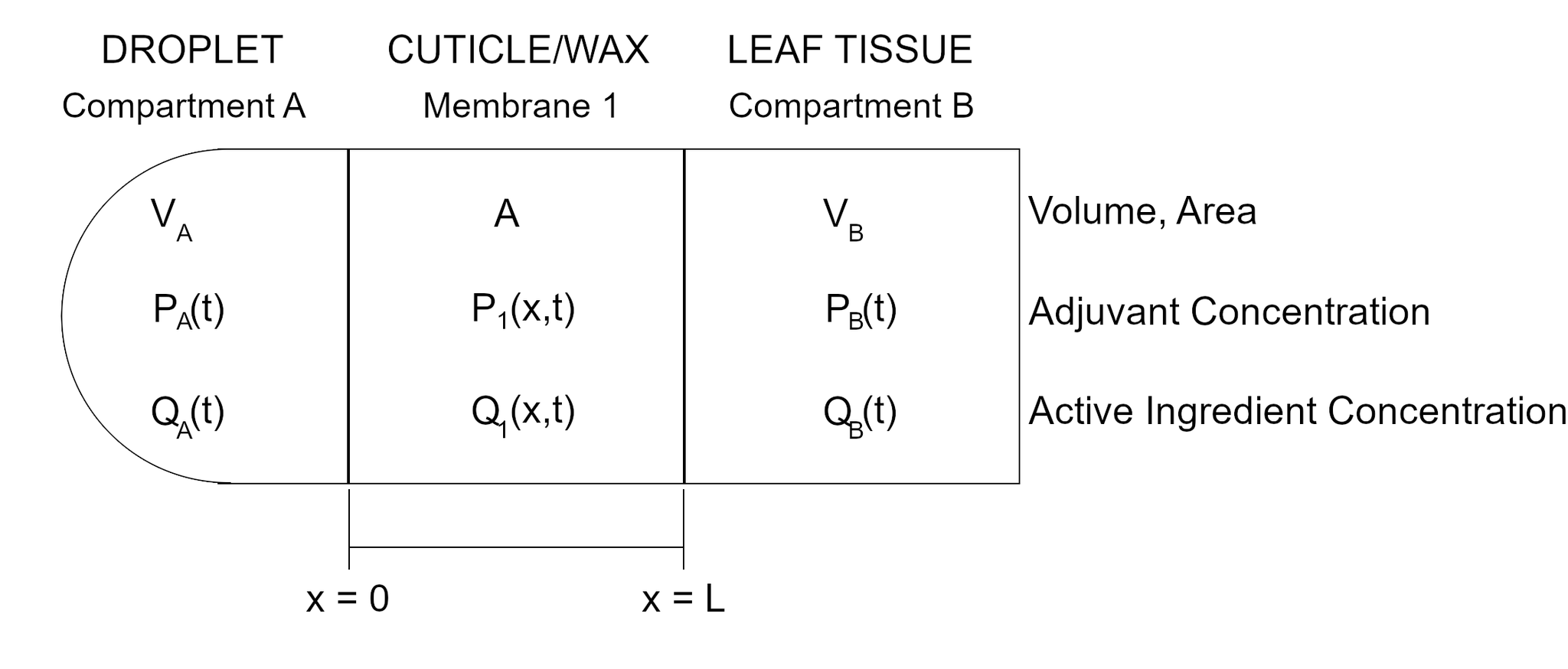}
		\centering
		\caption{Schematic diagram of the model system at time $t$ where state variables are defined in \cref{tab: nomenclature}.}
		\label{fig: hybrid_model}
	\end{figure}

        As for leaf anatomy, we use a simplified description of the pesticide formulation using two components - AI and AJ. The simplification lies in our decision to combine all compounds in the formulation which act to enhance uptake of AI across the leaf cuticle into a single variable AJ. \Cref{fig: hybrid_model} summarises our model structure.
			
        The leaf cuticle acts as a barrier for the plant and so, as with other biological barriers such as human skin, flow of solvents across the barrier is naturally slow \cite{SchreiberLukas2009WaSP,LaneMajellaE2013Spe}. This property determines our use of a spatially explicit system within the leaf cuticle. By contrast, in the droplet and the leaf tissue, we assume that both AI and AJ are well-mixed and consequently we describe them using ordinary differential rate equations. Furthermore, we assume that the droplet and leaf tissue have similar aqueous properties (relevant in our parameter estimation), that there is no crystallisation or photodegredation in the droplet and no metabolism within the leaf tissue and hence loss of AI and AJ from our system only occurs from the leaf tissue as the compounds move elsewhere in the leaf structure. 
			
        We describe flow between the three model compartments with partitioning \cite{HistoryPhysiological} which corresponds to the rate of flow being proportional to a weighted difference between current concentrations in the neighbouring compartments. The weightings are determined by the physico-chemical properties of the compound such that, flow between compartments ceases when the ratio of concentrations in each compartment equals to the partitioning coefficient. 
			
        As mentioned above, in the leaf cuticle we use a one-dimensional  spatio-temporal model. %Following Occam's Razor, 
        We model flow of the AJ in the cuticle using a simple diffusion process with a constant diffusion coefficient. We also use a diffusion process to describe the flow of AI across the leaf cuticle, but in this case we assume that the diffusion coefficient depends on the local concentration of the AJ. In particular, we choose a saturating function for the diffusion coefficient as described below. 

        Using the state variables presented in \cref{fig: hybrid_model} and defined in \cref{tab: nomenclature} and noting that 
	\begin{align*}
	    M_1(x,t)=A P_1(x,t) 
	\end{align*}
	our model system for the adjuvant is given as:
	\begin{subequations}
		\begin{align}
			\label{eq: PA_ode}
			\frac{d (V_A P_A)(t)}{dt} &= -\lambda_A A \left[ P_A(t) - \kappa_{A,1} \frac{M_1(0,t)}{A} \right], \\
			\label{eq: M1_pde}
			\frac{\partial M_1(x,t)}{\partial t} &= D_{P} \frac{\partial^2 M_1(x,t)}{\partial x^2}, \\
			\label{eq: PB_ode}
			\frac{d (V_B P_B)(t)}{dt} &= \lambda_B A \left[ \kappa_{B,1} \frac{M_1(L,t)}{A} - P_B(t) \right] - \beta V_B P_B(t).
		\end{align}
	\end{subequations}
        where all model parameters, defined in \cref{tab: nomenclature} are positive. To fully specify the system, and to ensure no loss of material across the membrane boundaries $x=0$ and $x=L$, we impose initial and boundary conditions:
	\begin{subequations}
		\begin{align}
			\label{eq: PA_ic}
			P_A(0) &= P_A^0 ,\\
			\label{eq: M1_ic}
			M_1(x,0) &= 0 ,\\
			\label{eq: PB_ic}
			P_B(0) &= 0 ,\\
			\label{eq: M1_bc_0}
			-D_{P} \frac{\partial M_1(0,t)}{\partial x} &= \lambda_A A \left[ P_A(t) - \kappa_{A,1} \frac{M_1(0,t)}{A} \right], \\
			\label{eq: M1_bc_L}
			-D_{P} \frac{\partial M_1(L,t)}{\partial x} &= \lambda_B A \left[ \kappa_{B,1} \frac{M_1(L,t)}{A} - P_B(t) \right]. 
		\end{align}
	\end{subequations}
			
        The model system for AI follows a similar structure except that the diffusion coefficient for the AI within the cuticle is assumed to be a saturating function of the local AJ concentration. This represents the facilitating role that the AJ plays in moving the AI into the leaf tissue. 

\begin{figure}[H]
		\includegraphics[width=7.5cm]{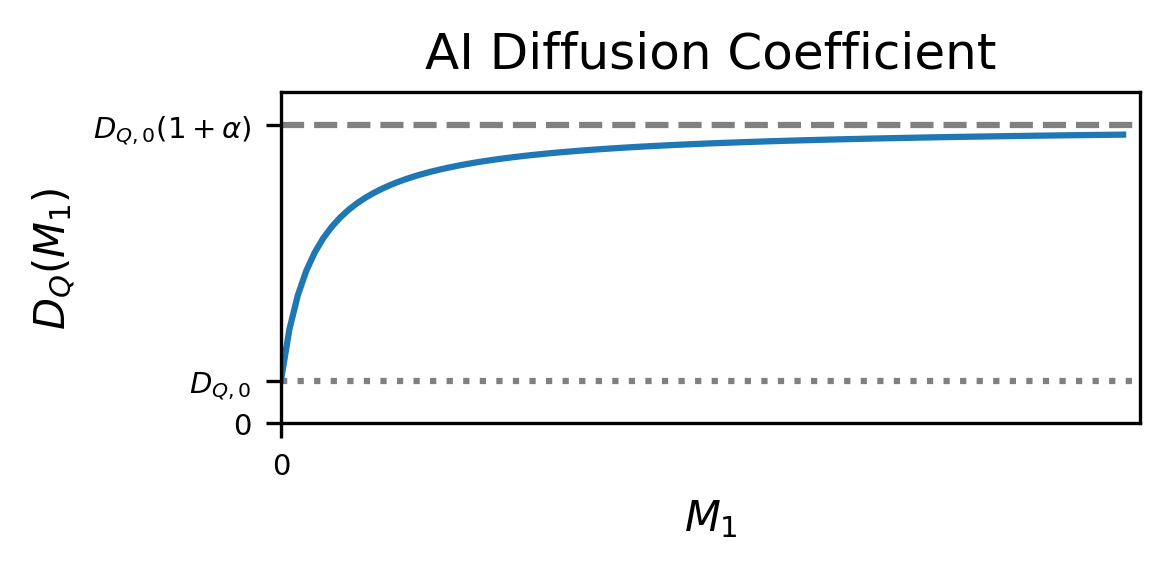}
		\centering
		\caption{Profile of AI diffusion coefficient as function of AJ concentration.}
		\label{fig: DQ}
	\end{figure}

        In the absence of any data, we assume that the concentration-dependent diffusion coefficient satisfies the parsimonious conditions:
	\begin{itemize}
		\item In the absence of AJ, the AI will still diffuse but more slowly.
		\item The effect of AJ is saturating such that at high AJ concentration, the impact on diffusion of AI is limited.
		\item The impact of AJ on diffusion of AI increases monotonically with AJ concentration.
	\end{itemize}
        An example of such a function which we use in our numerical simulations of the system is given by the algebraic expression 
	\begin{equation}
		\label{eq: DQ}
		D_Q(M_1) = D_{Q,0} \left[ 1 + \frac{\alpha M_1}{\sigma + M_1} \right],
	\end{equation}
        where $D_{Q,0}$, $\alpha$ and $\sigma$ are positive constants as defined in \cref{tab: nomenclature}. The shape of this function is shown in \cref{fig: DQ}.
	
	As above, we set 
	\begin{align*}
		N_1(x,t) &= A Q_1(x,t). 
	\end{align*}
        where $N_1$ denotes the amount of AI at a depth $x$ in the leaf cuticle and $Q_1$ is the concentration at that location. Using $N_1$ as the state variable within the model equation system we have:
	\begin{subequations}
		\begin{align}
			\label{eq: QA_ode}
			\frac{d (V_A Q_A)(t)}{dt} &= -\mu_A A \left[ Q_A(t) - K_{A,1} \frac{N_1(0,t)}{A} \right], \\
			\label{eq: N1_pde}
			\frac{\partial N_1(x,t)}{\partial t} &= \frac{\partial}{\partial x} \left[ D_{Q}(M_1) \frac{\partial N_1(x,t)}{\partial x} \right], \\
			\label{eq: QB_ode}
			\frac{d (V_B Q_B)(t)}{dt} &= \mu_B A \left[ K_{B,1} \frac{N_1(L,t)}{A} - Q_B(t) \right] - \eta V_B Q_B(t),	
		\end{align}
	\end{subequations}
        where $D_{Q}(M_1)$ is as described above. The system is fully specified with initial and boundary conditions: 
	\begin{subequations}
	    \begin{align}
			\label{eq: QA_ic}
			Q_A(0) &= Q_A^0 ,\\
			\label{eq: N1_ic}
			N_1(x,0) &= 0 ,\\
			\label{eq: QB_ic}
			Q_B(0) &= 0 ,\\
			\label{eq: N1_bc_0}
			-D_{Q}(M_1) \frac{\partial N_1(0,t)}{\partial x} &= \mu_A A \left[ Q_A(t) - K_{A,1} \frac{N_1(0,t)}{A} \right], \\
			\label{eq: N1_bc_L}
			-D_{Q}(M_1) \frac{\partial N_1(L,t)}{\partial x} &= \mu_B A \left[ K_{B,1} \frac{N_1(L,t)}{A} - Q_B(t) \right].
		\end{align}
	\end{subequations}

	\begin{table}[H]
    	\centering
		\caption{Summary of the model variables with descriptions and units. We use minutes ($min$) and micrometer ($\mu m$) for time and spatial parameters respectively, and percentage ($\%$) of the total amount for the amounts in the compartments.}
		\begin{tabular}{|c|c|c|}
			\hline
			Variable & Description & Unit \\ \hline
			$r$ & droplet radius & $\mu m$ \\ \hline
			$V_A$ & droplet volume & ${\mu m}^3$ \\ \hline
			$A$ & surface contact area & ${\mu m}^2$ \\ \hline
			$L$ & wax/cuticle length & $\mu m$ \\ \hline
			$L_B$ & leaf tissue length & $\mu m$ \\ \hline
			$V_B$ & leaf tissue volume & ${\mu m}^3$ \\ \hline
			$x$ &  distance in cuticle from droplet surface  &  $\mu m$ \\ \hline
			$t$ &  time variable  & $min$ \\ \hline
			$P_A(t)$ &  AJ droplet concentration at time $t$ & $\%\ {\mu m}^{-3}$ \\ \hline
			$P_1(x,t)$ &  AJ cuticle concentration at time $t$  & $\%\ {\mu m}^{-3}$ \\ \hline
			$P_B(t)$ &  AJ leaf tissue concentration at time $t$  & $\%\ {\mu m}^{-3}$ \\ \hline
			$M_1(x,t)$ &  AJ amount at depth $x$ at time $t$  & $\%\ {\mu m}^{-2}$ \\ \hline
			$Q_A(t)$ &  AI droplet concentration at time $t$ & $\%\ {\mu m}^{-3}$ \\ \hline
			$Q_1(x,t)$ &  AI cuticle concentration at time $t$  & $\%\ {\mu m}^{-3}$ \\ \hline
			$Q_B(t)$ &  AI leaf tissue concentration at time $t$  & $\%\ {\mu m}^{-3}$ \\ \hline
			$N_1(x,t)$ &  AI amount at depth $x$ at time $t$  & $\%\ {\mu m}^{-2}$ \\ \hline
			$D_P$ &  AJ diffusion coefficient in the cuticle  & ${\mu m}^2\ {min}^{-1}$ \\ \hline
			$D_Q$ &  AI diffusion coefficient in the cuticle  & ${\mu m}^2\ {min}^{-1}$ \\ \hline
			$D_Q^0$ &  AI baseline diffusion  & ${\mu m}^2\ {min}^{-1}$ \\ \hline
			$\alpha$ & saturation constant & $-$ \\ \hline
			$\sigma$ & positive constant & $\%\ {\mu m}^{-2}$ \\ \hline
			$\kappa_{A,1}$ & \makecell{AJ partition coefficient between droplet and cuticle, \\ note that $\kappa_{A,1}=\frac{1}{\kappa_{1,A}}$ }  &  $-$ \\ \hline
			$\kappa_{B,1}$ & \makecell{AJ partition coefficient between leaf tissue and cuticle, \\ note that $\kappa_{B,1}=\frac{1}{\kappa_{1,B}}$ } & $-$ \\ \hline
			$K_{A,1}$ & \makecell{AI partition coefficient between droplet and cuticle, \\ note that $K_{A,1}=\frac{1}{K_{1,A}}$ }  &  $-$ \\ \hline
			$K_{B,1}$ & \makecell{AI partition coefficient between leaf tissue and cuticle, \\ note that $K_{B,1}=\frac{1}{K_{1,B}}$ } & $-$ \\ \hline
			$\lambda_A$ & AJ speed of flow across droplet-cuticle boundary & ${\mu m}\ {min}^{-1}$ \\ \hline
			$\lambda_B$ & AJ speed of flow across cuticle-leaf tissue boundary & ${\mu m}\ {min}^{-1}$  \\ \hline
			$\mu_A$ & AI speed of flow across droplet-cuticle boundary & ${\mu m}\ {min}^{-1}$ \\ \hline
			$\mu_B$ & AI speed of flow across cuticle-leaf tissue boundary & ${\mu m}\ {min}^{-1}$  \\ \hline
			$\beta$ & AJ transfer rate from leaf tissue to rest of plant &  ${min}^{-1}$  \\ \hline
			$\eta$ & AI transfer rate from leaf tissue to rest of plant &  ${min}^{-1}$  \\ \hline
		\end{tabular}
		\label{tab: nomenclature}
	\end{table}
	
	\subsection{Experimental data}
		\label{subsec: Experimental Data}
\begin{figure}[H]
			\centering
			\begin{subfigure}{0.49\textwidth}
				\includegraphics[width=\textwidth]{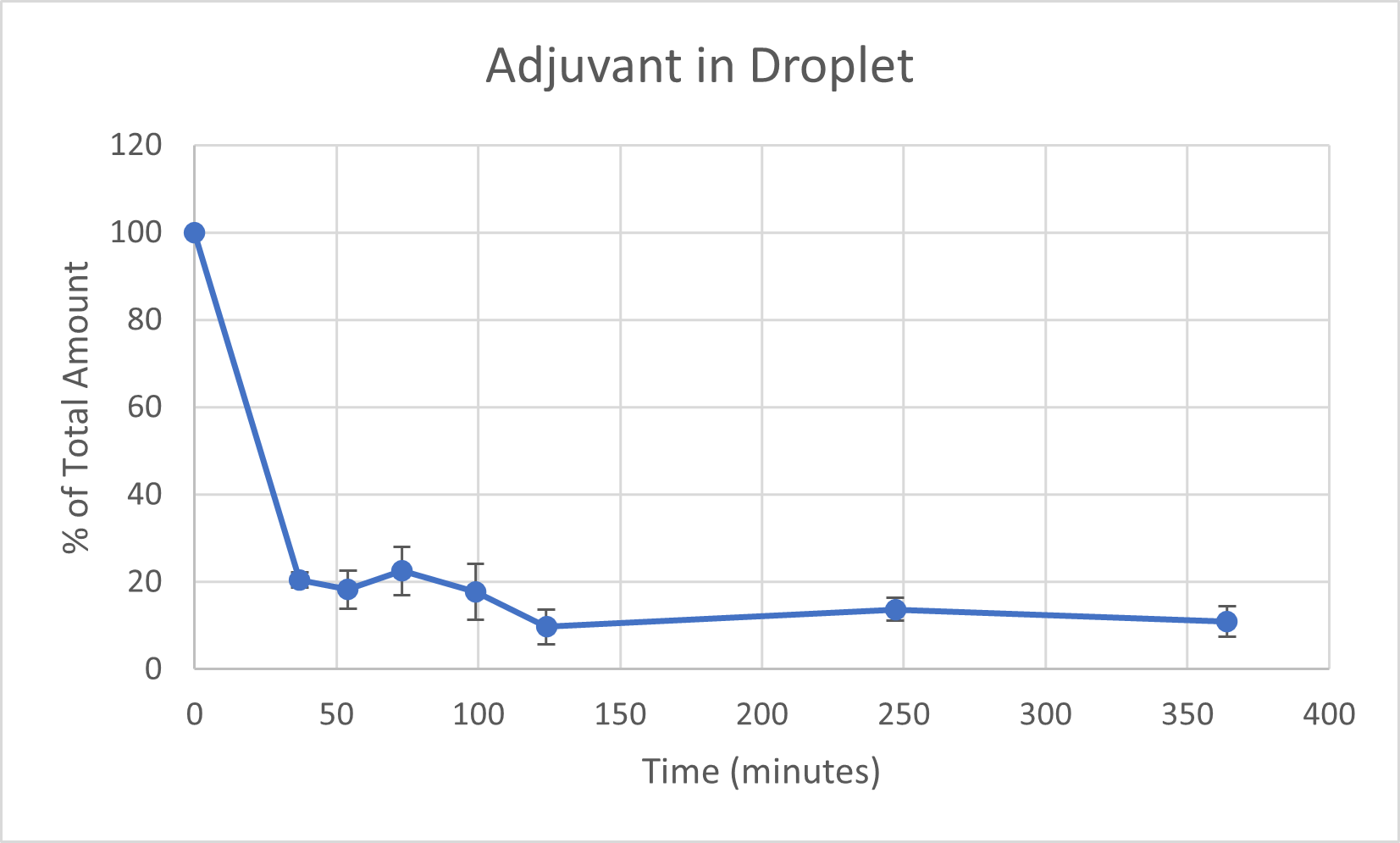}
				\caption{} 
				\label{fig: Adjuvant_droplet}
			\end{subfigure}
			\hfill 
			\begin{subfigure}{0.49\textwidth}
				\includegraphics[width=\textwidth]{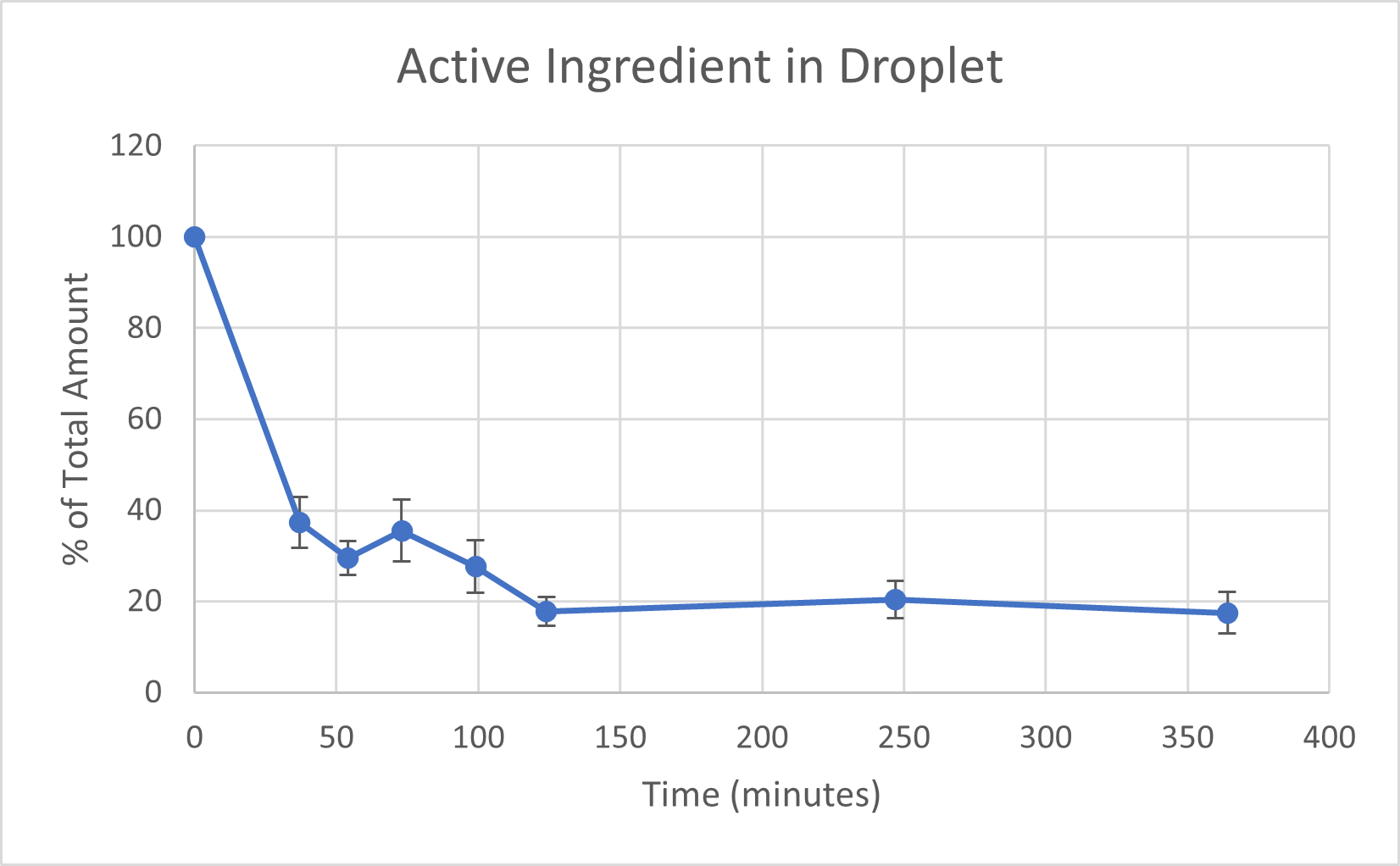}
				\caption{} 
				\label{fig: ActiveIngredient_droplet}
			\end{subfigure}
			\hfill 
			\begin{subfigure}{0.49\textwidth}
				\includegraphics[width=\textwidth]{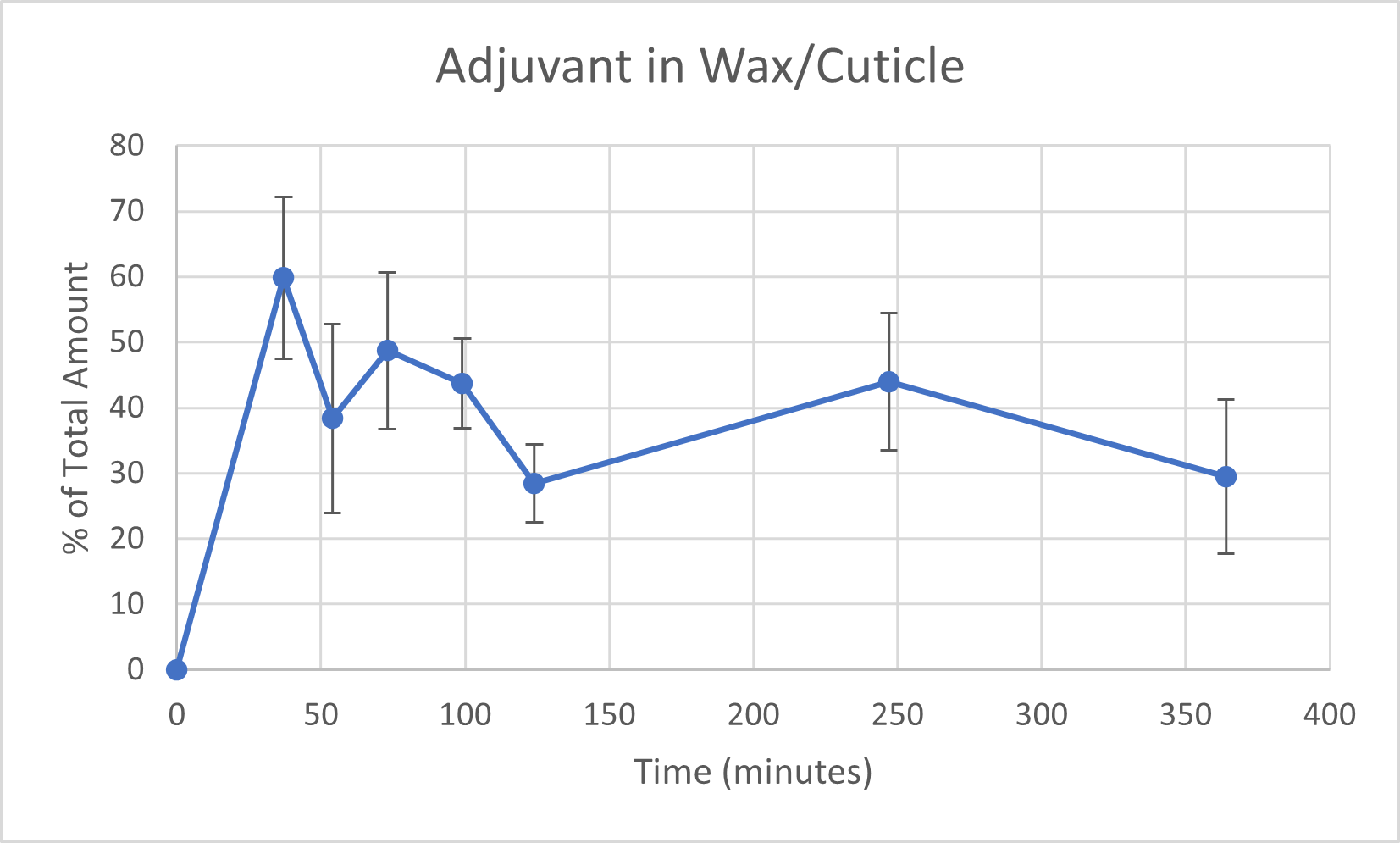}
				\caption{} 
				\label{fig: Adjuvant_cuticle}
			\end{subfigure}
			\hfill
			\begin{subfigure}{0.49\textwidth}
				\includegraphics[width=\textwidth]{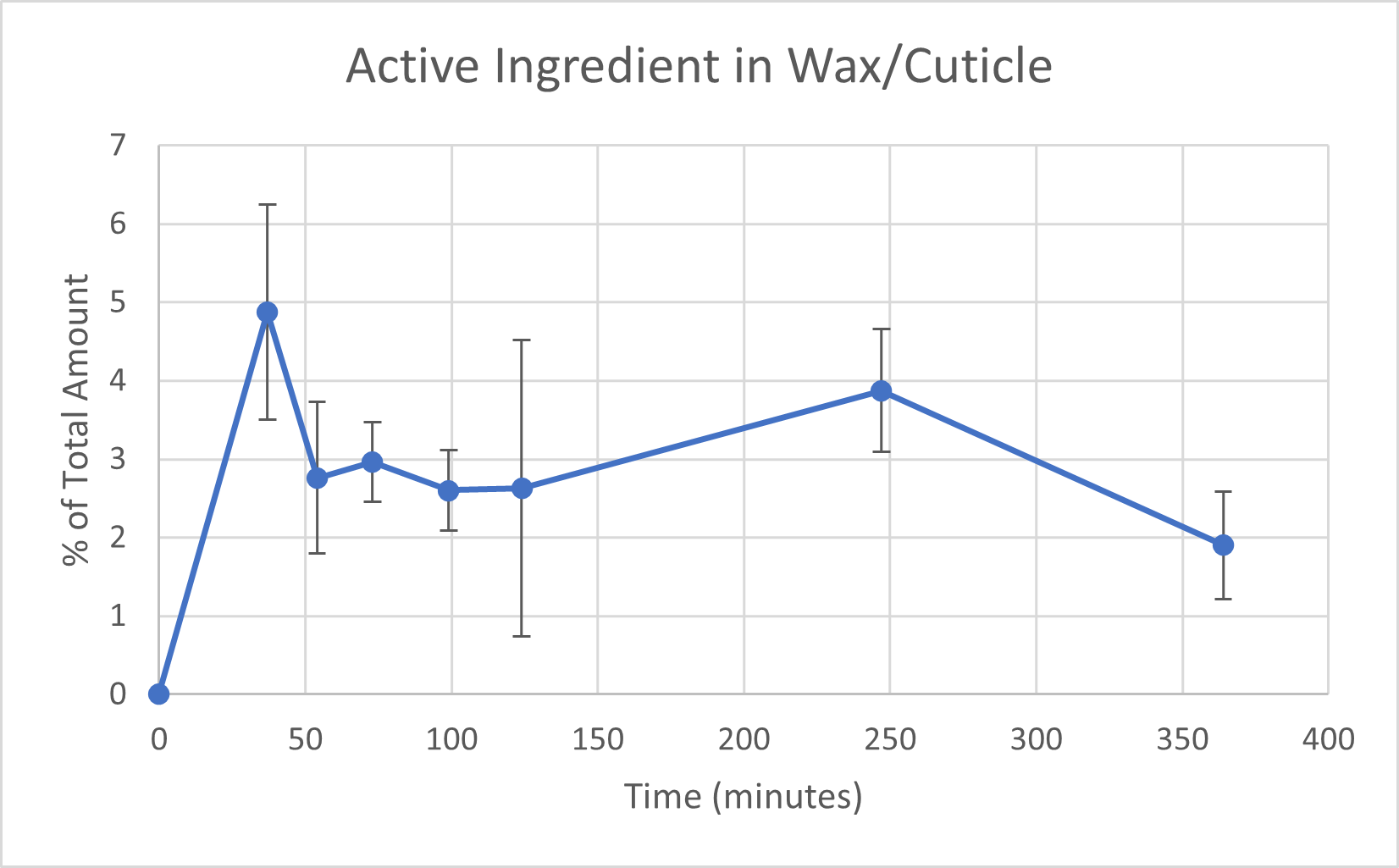}
				\caption{} 
				\label{fig: ActiveIngredient_cuticle}
			\end{subfigure}
			\hfill 
			\begin{subfigure}{0.49\textwidth}
				\includegraphics[width=\textwidth]{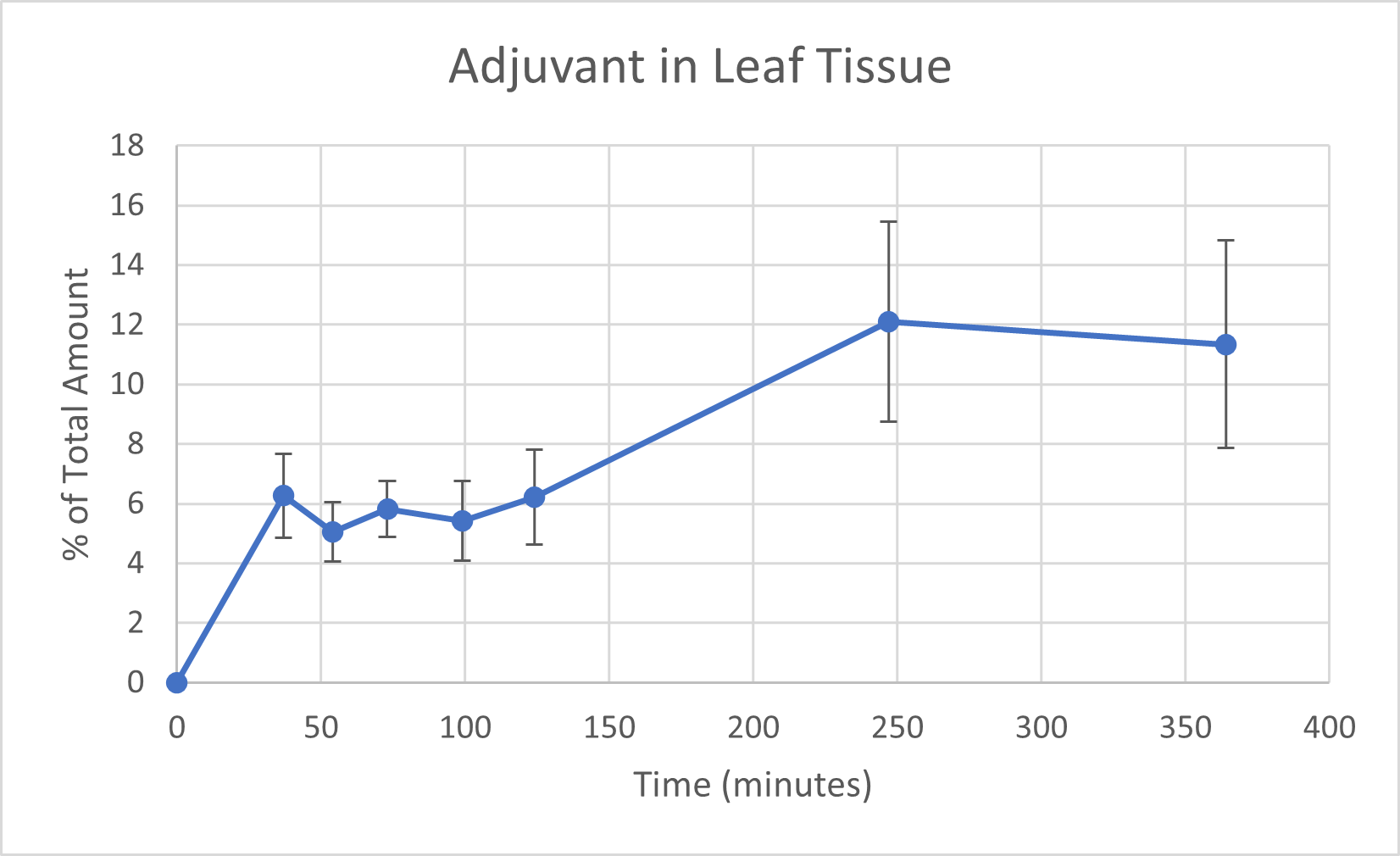}
				\caption{} 
				\label{fig: Adjuvant_leaftissue}
			\end{subfigure}
			\hfill 
			\begin{subfigure}{0.49\textwidth}
				\includegraphics[width=\textwidth]{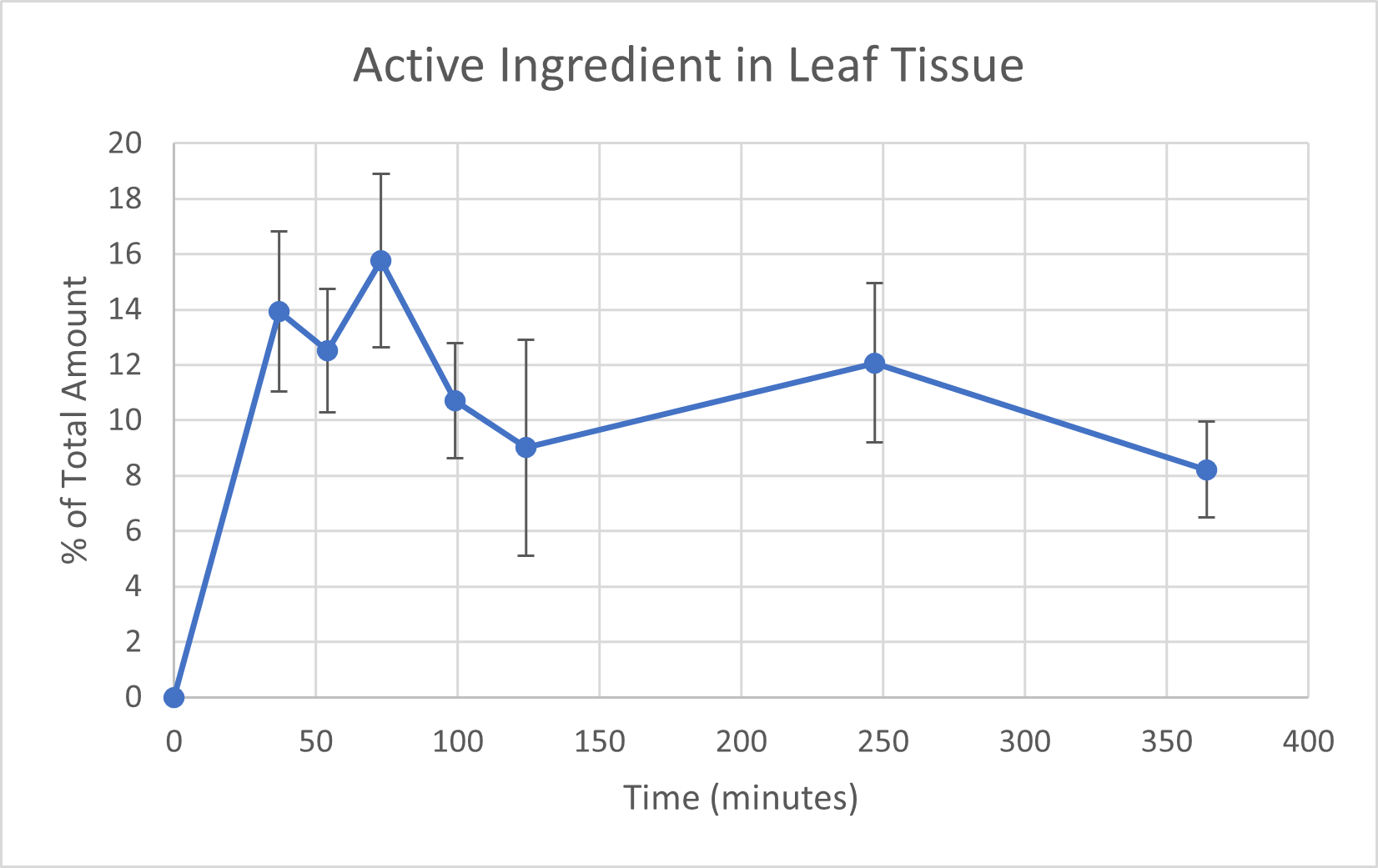}
				\caption{} 
				\label{fig: /ActiveIngredient_leaftissue}
			\end{subfigure}
			\hfill 
			\begin{subfigure}{0.49\textwidth}
				\includegraphics[width=\textwidth]{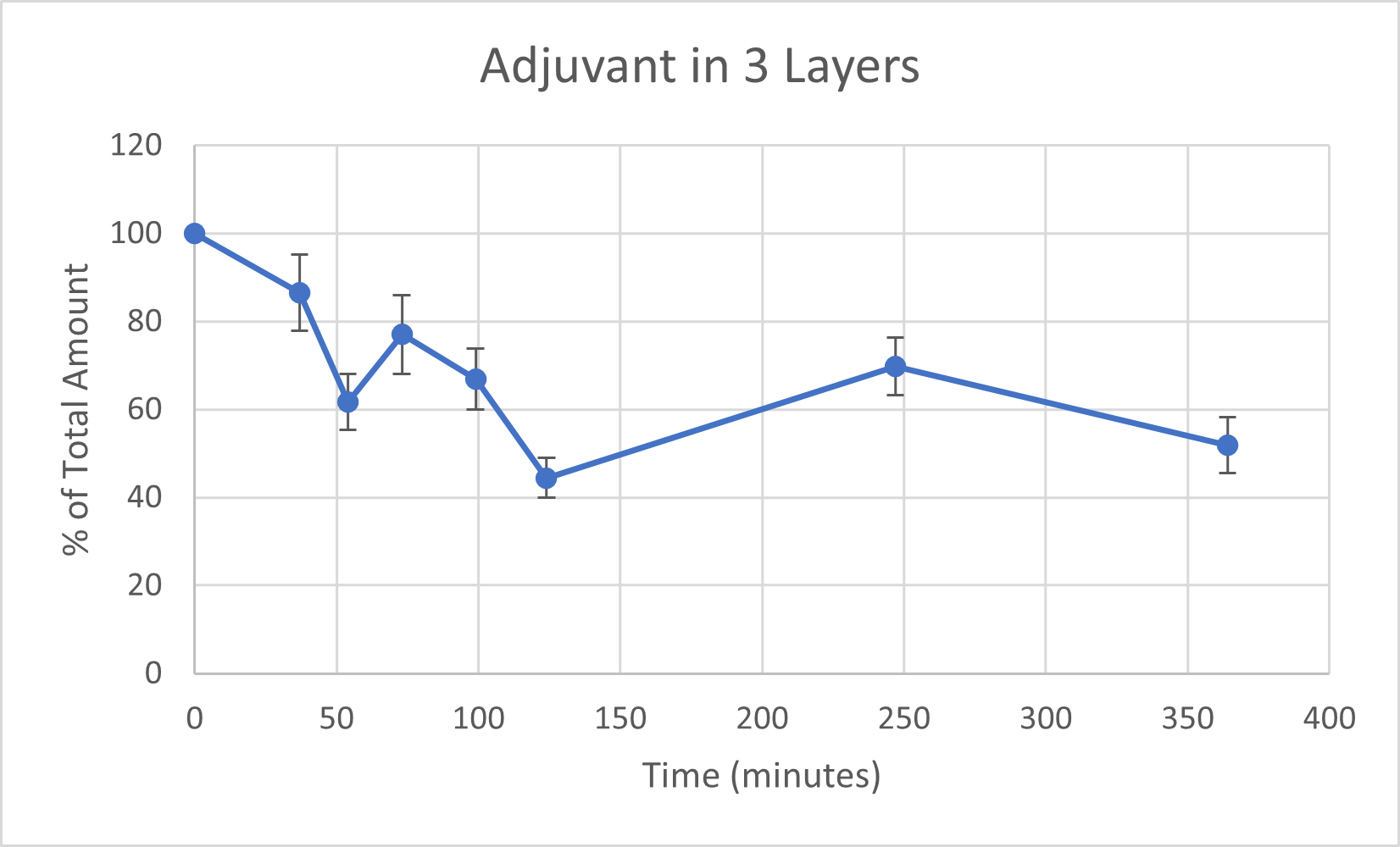}
				\caption{} 
				\label{fig: Adjuvant_layers}
			\end{subfigure}
			\hfill 
			\begin{subfigure}{0.49\textwidth}
				\includegraphics[width=\textwidth]{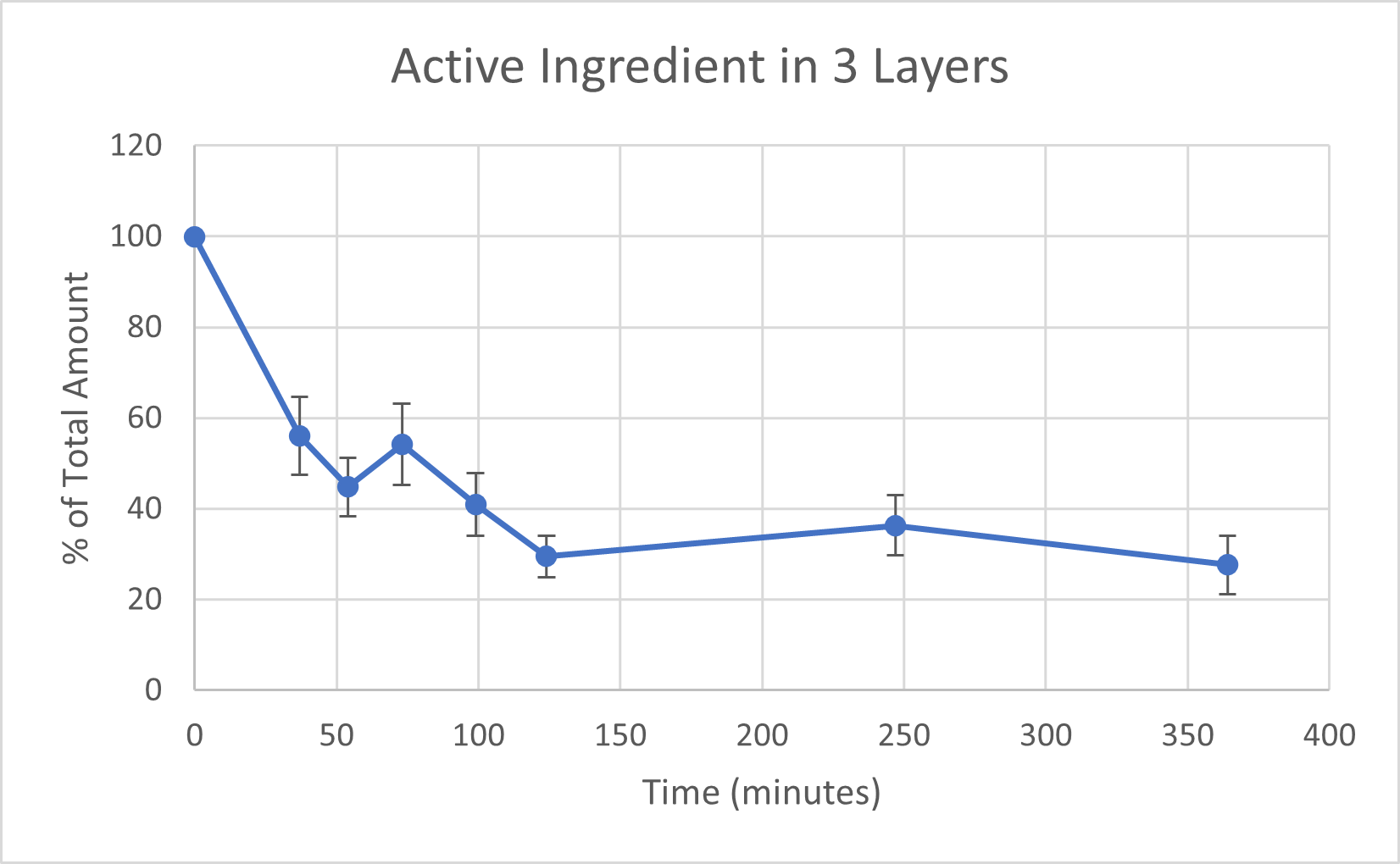}
				\caption{} 
				\label{fig: ActiveIngredient_layers}
			\end{subfigure}
                \caption{Means of the percentage of total amount initially available over time in minutes with 95\% CIs in each leaf layer for AJ ((a),(c),(e),(g)) and AI ((b),(d),(f),(h)).}
			\label{fig: ExperimentalData1}
		\end{figure}

            We use data collected by Syngenta \cite{Syngenta2022} and shown as a time series in \cref{fig: ExperimentalData1} to help parameterise and validate our model. At each data point we show the mean and 95\% confidence interval (CI) for 8 independent replicates separately for the AJ and AI. The data was obtained over a period of approximately 6 hours under conditions of no evaporation. Note that the experimental data records the percentage of the total amount in each compartment at the sample time. %which was the motivation for us to consider amounts rather than concentrations in the model formulation.			

            Droplet and leaf measurements corresponding to the experimental data have also been estimated at Syngenta; these estimates are provided in \cref{tab: space_time_parameters}.		
		\begin{table}[H]
			\centering
			\caption{Space and time parameter estimates based on experimental data provided Syngenta.}
			\begin{subtable}[h]{\textwidth}
				\centering
				\caption{\textit{Syngenta} data}
				\begin{tabular}{|c|c|c|c|}
					\hline
					Parameter & Description & Value & Unit \\ \hline
					$r$ & droplet radius & 30 & $\mu m$ \\ \hline
					$L$ & wax/cuticle length & 4 & $\mu m$ \\ \hline
					$L_B$ & leaf tissue thickness & 1000 & $\mu m$ \\ \hline
					$t_{end}$ & Experiment final time & 364 & $min$ \\ \hline
				\end{tabular}
				\label{tab: space_time_parameters_exp}
			\end{subtable}
			\vfill
			\vspace{.3in}
			\begin{subtable}[h]{\textwidth}
				\centering
				\caption{Calculated data}
				\begin{tabular}{|c|c|c|c|}
					\hline
					Parameter & Description & Value & Unit \\ \hline
					$V_A$ & droplet volume & $5.65\times10^{4}$ & ${\mu m}^3$ \\ \hline
					$A$ & surface contact area & $2.83\times10^3$ & ${\mu m}^2$ \\ \hline
					$V_B$ & leaf tissue volume & $2.83\times10^6$ & ${\mu m}^3$ \\ \hline
				\end{tabular}
				\label{tab: space_time_parameters_cal}
			\end{subtable}
			\label{tab: space_time_parameters}
		\end{table}

	\subsection{Parameter estimation}
		\label{subsec: Parameter Estimation}

            \subsubsection*{Diffusion coefficients ($D_P$, $D_{Q,0}$)}
                Since the volume of leaf tissue is significantly larger than either that of the droplet or the cuticle, we simplify our system to assume that the leaf tissue acts as a sink. In this case, the time $t_{lag}$ \cite{Edelstein-KeshetLeah2005Mmib} to traverse the cuticle of length $L$ is given by
                \begin{equation}
			    %\label{eq: diffusion_eq_01}
			    t_{lag} = \frac{L^2}{2 D},
			    \label{eq: diffusion_eq_02}
			\end{equation}
                where $D$ is diffusion coefficient.
                
                From \cref{fig: ExperimentalData1} we note that both AJ and AI have reached the leaf tissue by the time point $t=37\ min$, hence $0<t_{lag}<37\ min$. With no further data available, we choose the interval $5\leq t_{lag} \leq 20\ min$ to estimate $D$; our estimates are presented in \cref{tab: diffusion_coefficients}.

		      \begin{table}[H]
			    \centering
			    \caption{Approximations for the diffusion coefficients using \cref{eq: diffusion_eq_02} and data from \cref{fig: ExperimentalData1}.}
		      	\begin{tabular}{|c|c|c|}
			     	\hline
		  		& Range & unit \\ \hline
		  		$t_{lag}$ & $(5,20)$ & $min$ \\ \hline
		  		$D$ & (0.4, 1.6) & ${\mu m}^2/min$ \\ \hline
		      	\end{tabular}
		      	\label{tab: diffusion_coefficients}
	       	\end{table}

            \subsubsection*{Partition coefficients ($\kappa_{i,j}$, $K_{i,j}$)} 
                We estimate the partition coefficients using the last 3 time points in the experimental data since the coefficients give the ratio of concentrations across a membrane when those concentrations are in equilibrium. We further simplify our estimate by assuming that partitioning between the cuticle and either the droplet or the leaf tissue is the same to reflect the context that both are water-rich environments \cite{CannyM.J2006LWCa,Oguchi2018} (in sharp contrast to the waxy cuticle layer). We have 
			\begin{equation}
				\begin{array}{ll}
					\kappa_{1,i} = \displaystyle \frac{P_1}{P_i}, ~&~ 
					K_{1,i} = \displaystyle \frac{Q_1}{Q_i},
				\end{array}
				\label{eq: partition_eq}
			\end{equation}
                for $i=A,B$ and note that our assumption on water-rich environments means that $\kappa_{1,A}=\kappa_{1,B}$, and $K_{1,A}=K_{1,B}$. Estimates are presented in \cref{tab: partition_coefficients}.
				
                \begin{table}[H]
				\centering
				\caption{Approximations for the partition coefficients using  \cref{eq: partition_eq} and experimental data from \cref{fig: ExperimentalData1}.}
				\begin{tabular}{|c|c|c|c|}
					\hline
					&  AJ & Range & unit \\ \hline
					$P_1$ & 3.00$\times 10^{-3}$ & (2.51$\times 10^{-3}$, 3.50$\times 10^{-3}$) & $\%\ {\mu m}^{-3}$ \\ \hline
					$P_A$ & 2.03$\times 10^{-4}$ & (1.71$\times 10^{-4}$, 2.35$\times 10^{-4}$) & $\%\ {\mu m}^{-3}$ \\ \hline
					$\kappa_{1,A}$ & 14.80 & (10.67, 20.48) & $-$ \\ \hline
					&  &  &  \\ \hline
					&  AI & Range & unit \\ \hline
					$Q_1$ & 2.48$\times 10^{-4}$ & (1.86$\times 10^{-4}$, 3.09$\times 10^{-4}$) & $\%\ {\mu m}^{-3}$ \\ \hline
					$Q_A$ & 3.28$\times 10^{-4}$ & (2.93$\times 10^{-4}$, 3.64$\times 10^{-4}$) & $\%\ {\mu m}^{-3}$ \\ \hline
					$K_{1,A}$ & 7.54$\times 10^{-1}$ & (5.11$\times 10^{-1}$, 1.06) & $-$ \\ \hline
				\end{tabular}
				\label{tab: partition_coefficients}
			\end{table}

            \subsubsection*{Speeds ($\lambda_A$, $\mu_A$)}
                Since there is initially no chemical in the leaf cuticle, we approximate the initial loss of AI and AJ from the droplet by
                \begin{equation*}
				\label{eq: PA_ode_lamdaA}
				\frac{d (V_A X_A)(t)}{dt} = -y_A A X_A(t),
			\end{equation*}
                where $X_A = P_A$ or $Q_A$ and $y_A = \lambda_A$ or $\mu_A$. We assume that this simple exponential decay occurs until at least the first recorded time point $t=37$ minutes and obtain the mean and 95\% CI for each parameter by using the data from \cref{fig: ExperimentalData1} in the expression
                \begin{equation}
				\label{eq: PA_lambda}
	           	y_A = - \frac{V_A \times \ln \left( \frac{X_A(t)}{{X_A}^0} \right)}{A\times t}.
			\end{equation}
                The relevant parameter estimates are given in \cref{tab: lambda_mu_A}.
                \begin{table}[H]
				\centering
				\caption{Approximations for $\lambda_A$, $\mu_A$ using \cref{eq: PA_lambda}, and data from \cref{fig: ExperimentalData1} and \cref{tab: space_time_parameters}.}
				\begin{tabular}{|c|c|c|c|}
					\hline
					&  AJ & Range & unit \\ \hline
					$t$ & 37 & -- & $min$ \\ \hline
					${P_A}^0$ & 1.77$\times10^{-3}$ & -- & $\%\ {\mu m}^{-3}$ \\ \hline
					$P_A(t)$ & 3.62$\times10^{-4}$ & (3.32$\times10^{-4}$, 3.91$\times10^{-4}$) & $\%\ {\mu m}^{-3}$ \\ \hline
					$\lambda_A$ & 8.58$\times10^{-1}$ & (8.16$\times10^{-1}$, 9.04$\times10^{-1}$) & $\mu m/min$ \\ \hline
					&  &  &  \\ \hline
					&  AI & Range & unit \\ \hline
					$t$ & 37 & -- & $min$ \\ \hline
					${Q_A}^0$ & 1.77$\times10^{-3}$ & -- & $\%\ {\mu m}^{-3}$ \\ \hline
					$Q_A(t)$ & 6.69$\times10^{-4}$ & (5.60$\times10^{-4}$, 7.59$\times10^{-4}$) & $\%\ {\mu m}^{-3}$ \\ \hline
					$\mu_A$ & 5.33$\times10^{-1}$ & (4.57$\times10^{-1}$, 6.22$\times10^{-1}$) & $\mu m/min$ \\ \hline
				\end{tabular}
				\label{tab: lambda_mu_A}
			\end{table}
				
            \subsubsection*{Speeds ($\lambda_B$, $\mu_B$)}
                Unlike with $\lambda_A$ and $\mu_A$, there is no obvious timescale argument that we can make to estimate the parameters $\lambda_B$ and $\mu_B$. Therefore we once again take the parsimonious approach and assume that the speed with which molecules cross into and out of the cuticle is the same. In terms of parameter estimates, this means that we take $\lambda_B=\lambda_A$ and $\mu_B=\mu_A.$

            \subsubsection*{Transfer rates ($\beta$, $\eta$)}
                The two parameters $\beta$ and $\eta$ denote the loss rates of AJ and AI from the system into the leaf transport system. Using data from \cref{fig: ExperimentalData1} we can estimate for each time point, the fraction of material that has left the model system. Noting that 
                \begin{equation}
				\label{eq: beta_00}
				\frac{d}{dt} \left[ V_A X_A(t) + \int_{0}^{L} Y_1(x,t)+ V_B X_B(t) \right] = - \nu V_B X_B(t),
			\end{equation}
                where $(X_A,Y_1,X_B) \in \{(P_A, M_1, P_B), (Q_A, N_1, Q_B)\}$ and $\nu \in \{\beta, \eta \}$ we use the data points from \cref{fig: ExperimentalData1} to estimate the rate of change in amounts between time measurements. This results in the mean and 95\% CI for the two parameters as shown in \cref{tab: beta_eta}.

                \begin{table}[H]
				\centering
				\caption{Approximations for $\beta$, $\eta$ using \cref{eq: beta_00} with $\frac{d}{dt} \left[ \cdot \right]$ denoting its left-hand side, and results from \cref{fig: ExperimentalData1}.}
				\begin{tabular}{|c|c|c|c|}
					\hline
					&  AJ & Range & unit \\ \hline
					$\frac{d}{dt} \left[ \cdot \right]$ & 8.92$\times10^{-2}$ & $-$ & $\%\ {min}^{-1}$  \\ \hline
					$P_B(t)$ & 2.31$\times10^{-6}$ & (1.94$\times10^{-6}$, 2.68$\times10^{-6}$) & $\%\ {\mu m}^{-3}$  \\ \hline
					$\beta$ & 1.37$\times10^{-2}$ & (1.18$\times10^{-2}$, 1.63$\times10^{-2}$) & ${min}^{-1}$ \\ \hline
					&   &   &   \\ \hline
					&  AI & Range & unit \\ \hline
					$\frac{d}{dt} \left[ \cdot \right]$ & 1.29$\times10^{-1}$ & $-$ & $\%\ {min}^{-1}$  \\ \hline
					$Q_B(t)$ & 3.63$\times10^{-6}$ & (3.16$\times10^{-6}$, 4.11$\times10^{-6}$) & $\%\ {\mu m}^{-3}$  \\ \hline
					$\eta$ & 1.26$\times10^{-2}$ & (1.11$\times10^{-2}$, 1.45$\times10^{-2}$) & ${min}^{-1}$ \\ \hline
				\end{tabular}
				\label{tab: beta_eta}
			\end{table}

            \subsubsection*{Parameters associated with concentration-dependent diffusion}
                We do not estimate the two model parameters $\alpha$ and $\sigma$; rather we use them to explore the impact of concentration-dependent diffusion on the transport of the AI into the plant via leaf cuticle. 
			
    \section{Results}
	\label{sec: Results}
        Our results are focused on numerical solutions of our model system detailed in \cref{sec: Model} and how they are affected by changes in the concentration-dependent diffusion. However, to set the scene, we begin by presenting steady state solutions assuming that there is no flow from the leaf tissue into the rest of the plant structure.

        \subsection*{Steady-state analysis}
            With $\beta=\eta=0$ the amount of mass in our model     system is conserved. In this case, we can find the steady state distribution of AJ and AI by setting each rate equation equal to zero (and ensuring that the total at each time point equals the initial amount in the droplet). In doing this we arrive at the following algebraic expressions in each compartment: \\

            \begin{itemize}
                \item {\bf AJ steady state distributions:}
                    \begin{equation}
                        P_A = \frac{\kappa_{A,1} V_A {P_A}^0 }{\kappa_{A,1} V_A + \kappa_{B,1} V_B + A L}, \quad M_1 = \frac{A P_A}{\kappa_{A,1}}, \quad P_B = \frac{\kappa_{B,1}}{\kappa_{A,1}} P_A.
                    \end{equation}
                \item {\bf AI steady state distributions:}
                    \begin{equation}
                        \label{eq:ststAI}
                        Q_A = \frac{K_{A,1} V_A {Q_A}^0 }{K_{A,1} V_A + K_{B,1} V_B + A L}, \quad N_1 = \frac{A Q_A}{K_{A,1}}, \quad Q_B = \frac{K_{B,1}}{K_{A,1}} Q_A.    
                    \end{equation}
            \end{itemize}

		\begin{figure}[H]
			\centering
			\begin{subfigure}{0.32\textwidth}
				\includegraphics[width=\textwidth]{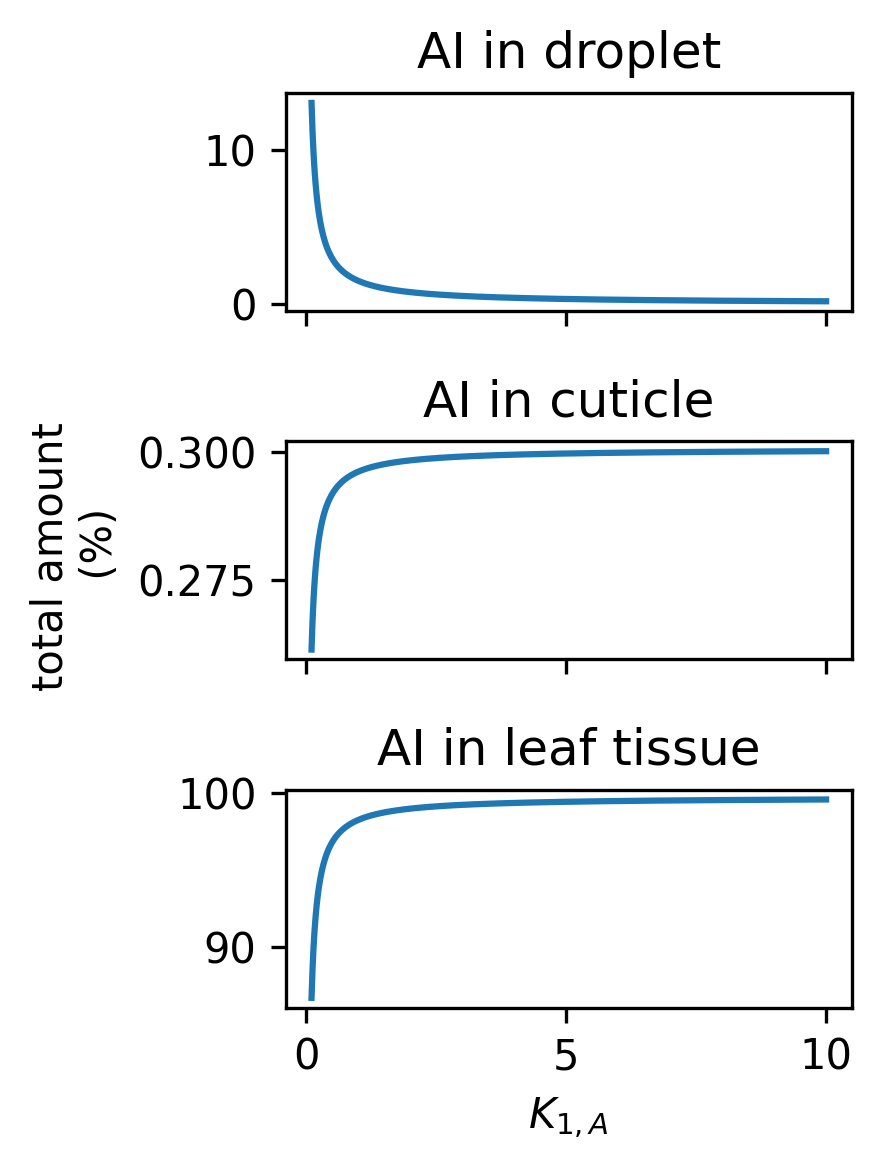}
				\caption{}
				\label{fig: TOTAL_per_SS_K_1A}
			\end{subfigure}
			\hfill
			\begin{subfigure}{0.32\textwidth}
				\includegraphics[width=\textwidth]{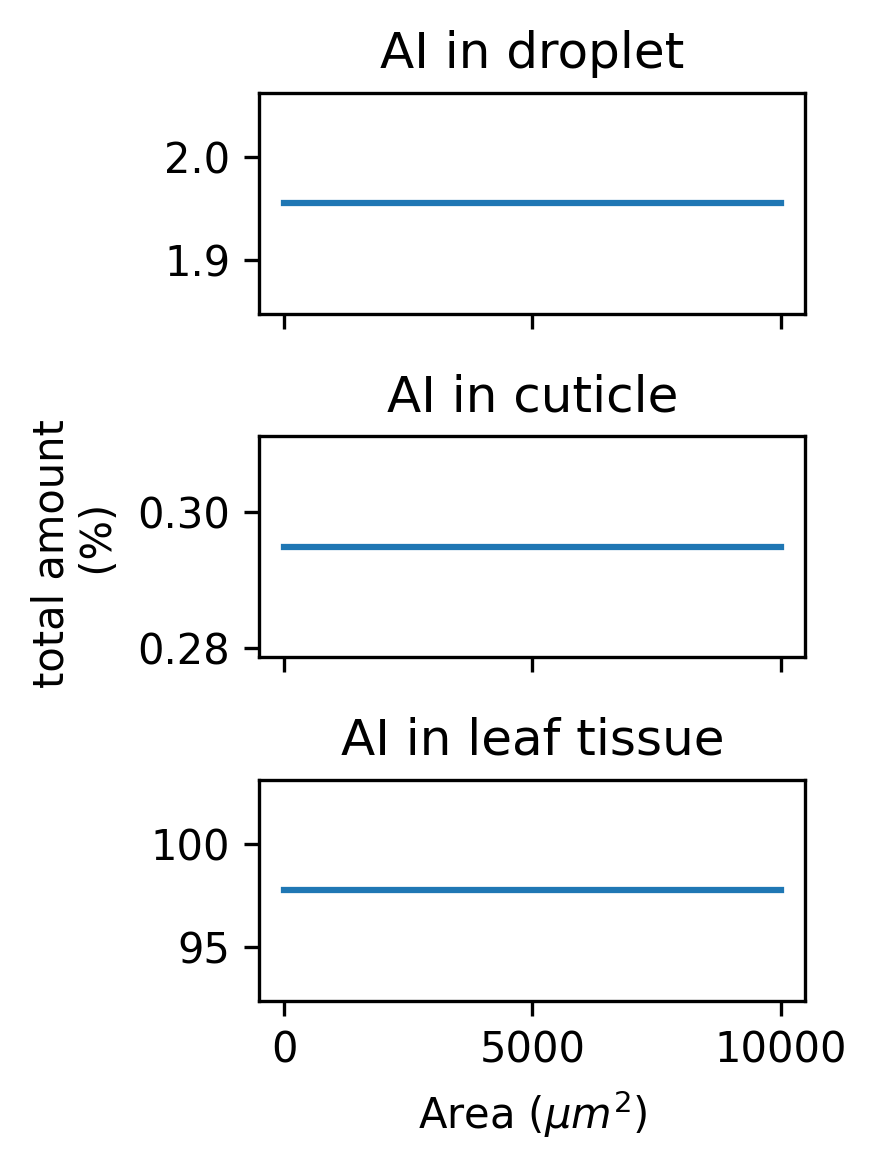}
				\caption{}
				\label{fig: TOTAL_per_SS_A}
			\end{subfigure}
			\hfill
			\begin{subfigure}{0.32\textwidth}
				\includegraphics[width=\textwidth]{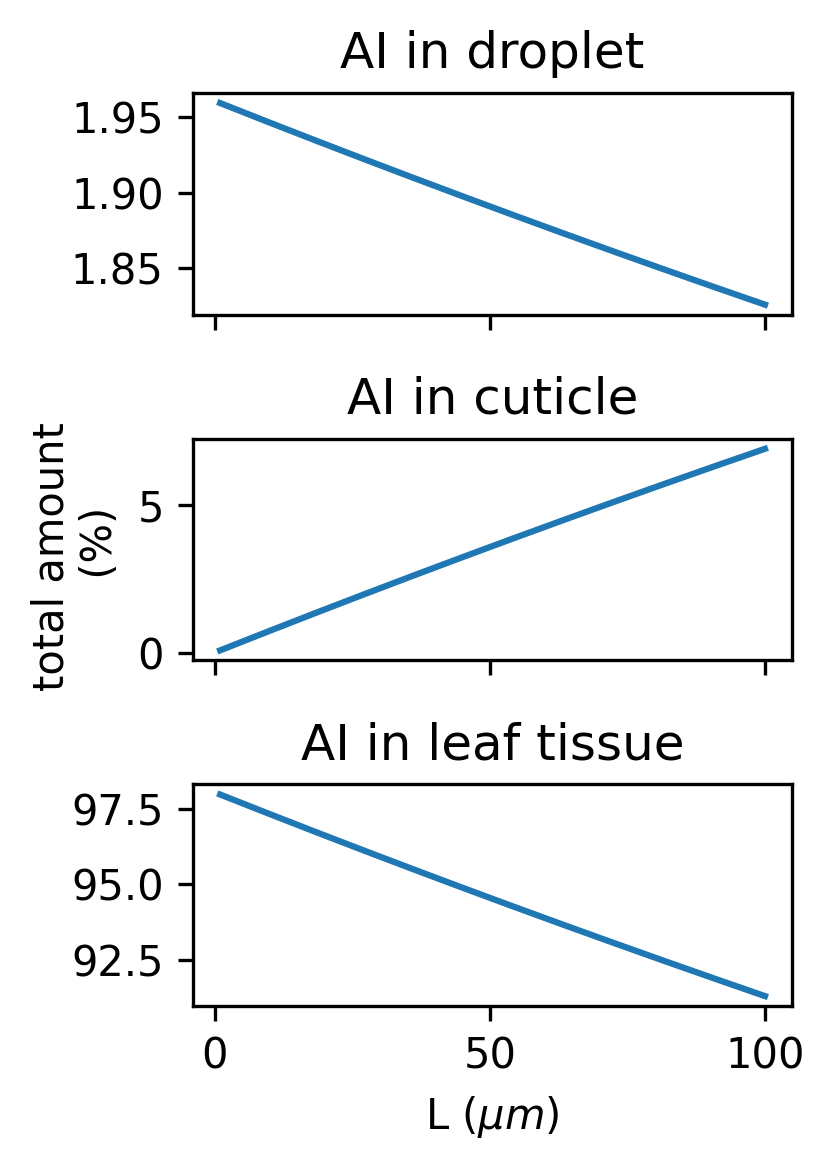}
				\caption{}
				\label{fig: TOTAL_per_SS_L}
			\end{subfigure}	
			\caption{Steady-state solution to the AI model system \cref{eq:ststAI}  using parameter estimates in \cref{subsec: Parameter Estimation} for a range of values for (a)$K_{1,A}$ , (b)$A$ , and (c)$L$.}
			\label{fig: TOTAL_per_SS}		
		\end{figure}	
            The plots shown in \cref{fig: TOTAL_per_SS} are intuitive: the greater the affinity for the chemical to be in a non-aqueous region, the greater the amount that will accumulate in the cuticle (a); the area of the sample region does not affect the percentage of total amount in that region since we assume homogeneity in the lateral direction (b); and as the length of cuticle increases, so does the amount of chemical in that compartment (c). One other point worth making here is that the impact of concentration-dependent concentration is not seen in the steady state so the impact of AJ will not be identified if the data is only provided for times after which the system is at steady state.

        \subsection{Time-dependent solutions}
            We solve our model system numerically using a Finite Element Method \cite{FEM3,FEM1} for the spatial discretisation of the cuticle, and explicit Finite Difference Methods \cite{FDM2,FDM1} for the time discretisation. The method was designed to conserve total mass by carefully discretising the boundary terms and rates of change, and to be stable even with very small diffusion coefficients. We use the Python package Firedrake \cite{RathgeberFlorian2017FAtF,MitchellL2019faaf} to facilitate this. The code can be accessed via Github \cite{DelosReyesGithubNoEvap}.
            
            %of the system. The major challenge arises from  our boundary conditions; we use \textit{Python} \cite{LangtangenHansPetter2009APoS,Python} with \textit{Spyder} \cite{Spyder} as a compiler, and a package called	\textit{Firedrake} \cite{RathgeberFlorian2017FAtF,MitchellL2019faaf} to facilitate this. The code can be accessed via Github \cite{DelosReyesGithubNoEvap}.

            We consider two cases: firstly when there is no impact of AJ on the diffusion of AI across the cuticle ($\alpha=0$) and then in the case where concentration-dependent diffusion does occur ($\alpha>0$).
			
		\subsubsection{Diffusion of AI independent of local AJ concentration}
\begin{figure}[H]
                    \centering
				\includegraphics[width=11cm]{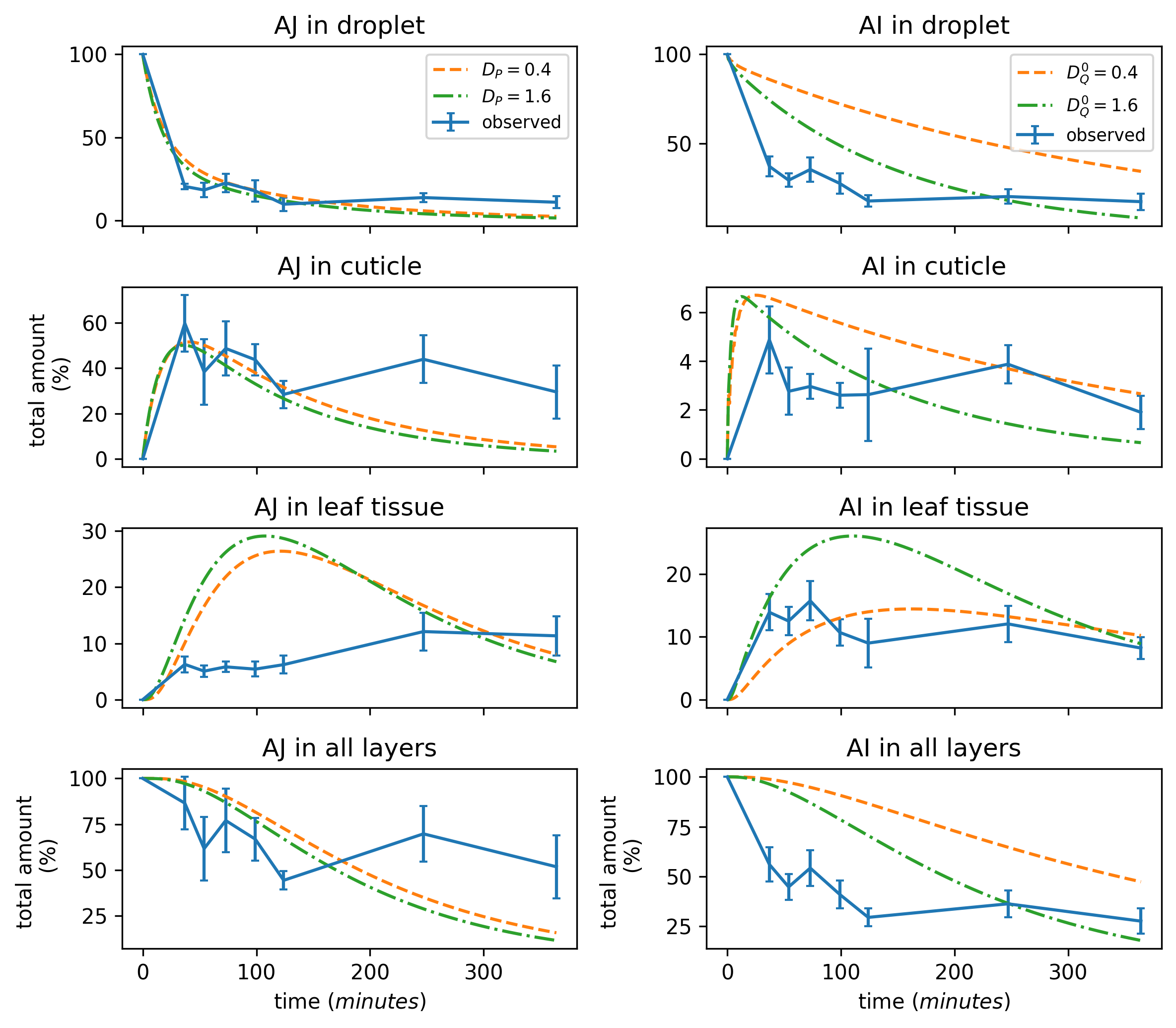}
				\caption{Predicted uptake of AJ and AI using parameter estimates from \cref{subsec: Parameter Estimation} with $D_P=0.4,1.6$, and $D_Q=0.4,1.6$ with $\alpha=0$.}
				\label{fig: Numerical_Solution_alpha0}
			\end{figure}

                With $\alpha=0$, we use the parameter estimates described in the previous \cref{subsec: Parameter Estimation} to solve our model system. The results are shown in \cref{fig: Numerical_Solution_alpha0} for the lower and upper bound on the constant diffusion coefficient. 

                We are encouraged that in several of the subplots, there is a reasonably good approximation to the data either across the full timespan (eg. AJ in the droplet) or during the initial phase (eg. AJ in the cuticle, AI in the leaf tissue with lower bound $D_Q^0$). There is an overestimation in the leaf tissue which is more apparent with AJ and AI with higher $D_Q^0$. We can observe in \cref{fig: Numerical_Solution_alpha0} that the range for AI diffusion coefficient is still rather slow to allow for more AI to move across the leaf layers and into the rest of the plant. This is not surprising since we have not yet accounted for the accelerating effects of the presence of the AJ.

            \subsubsection{Concentration dependent AI diffusion: $\alpha>0$, $\sigma>0$}
\begin{figure}[H]
				\centering
				\includegraphics[width=10cm]{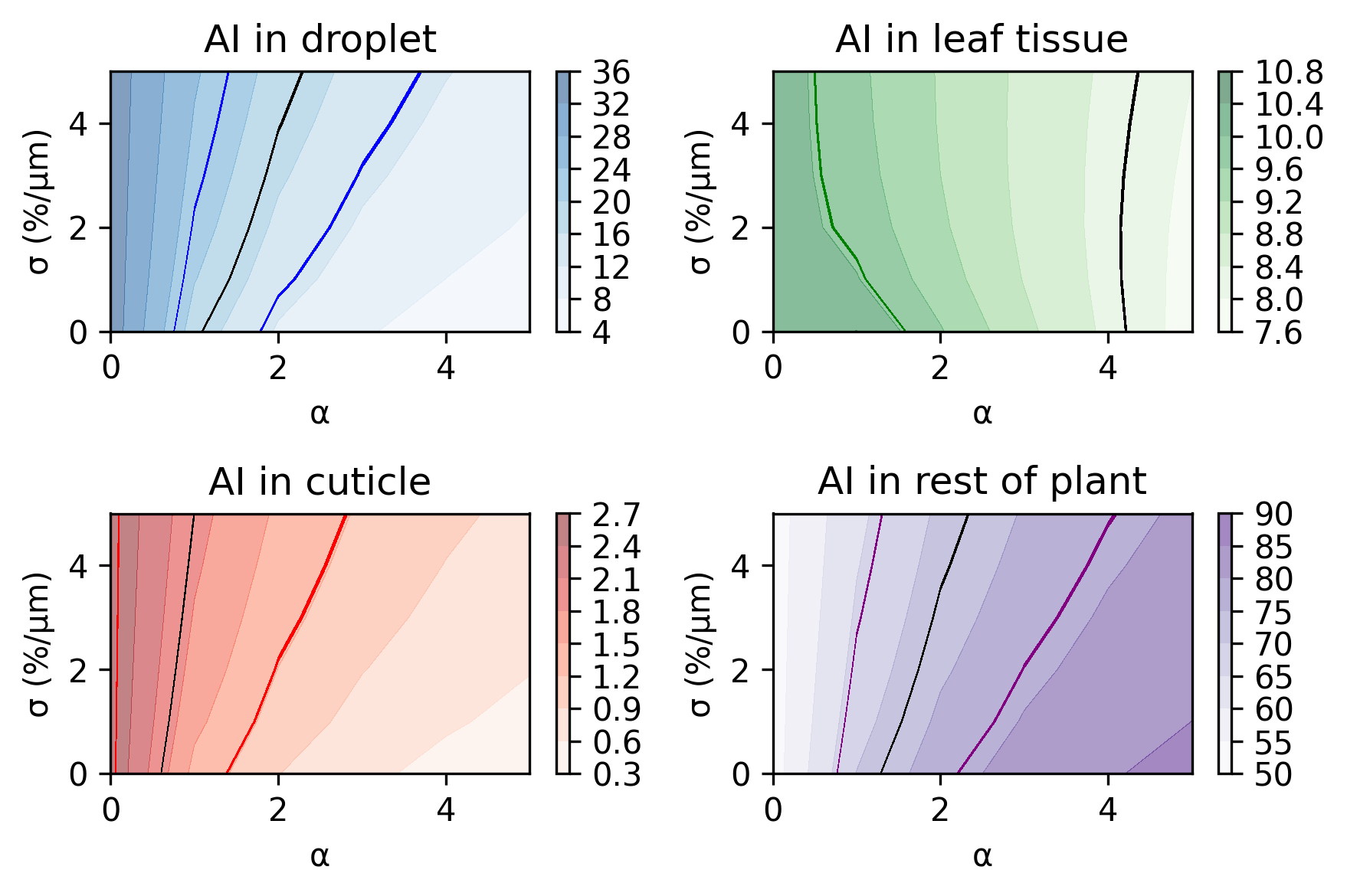}
				\caption{Predicted percentage of the total amount of the AI in droplet, leaf, and the rest of the plant which is the amount that has been removed from the system. We use parameter estimates in \cref{subsec: Parameter Estimation} with $D_P,D_Q^0=0.4$ for a range of values for $\alpha$ and $\sigma$ at time $t=364\ min$. The solid curves represent the means (black) and 95\% CIs (coloured) of the total amounts from the experimental data.}
				\label{fig: Numerical_Solution_alpha_sigma_DQ00.4}
			\end{figure}

\begin{figure}[H]
                    \centering
        		\includegraphics[width=10cm]{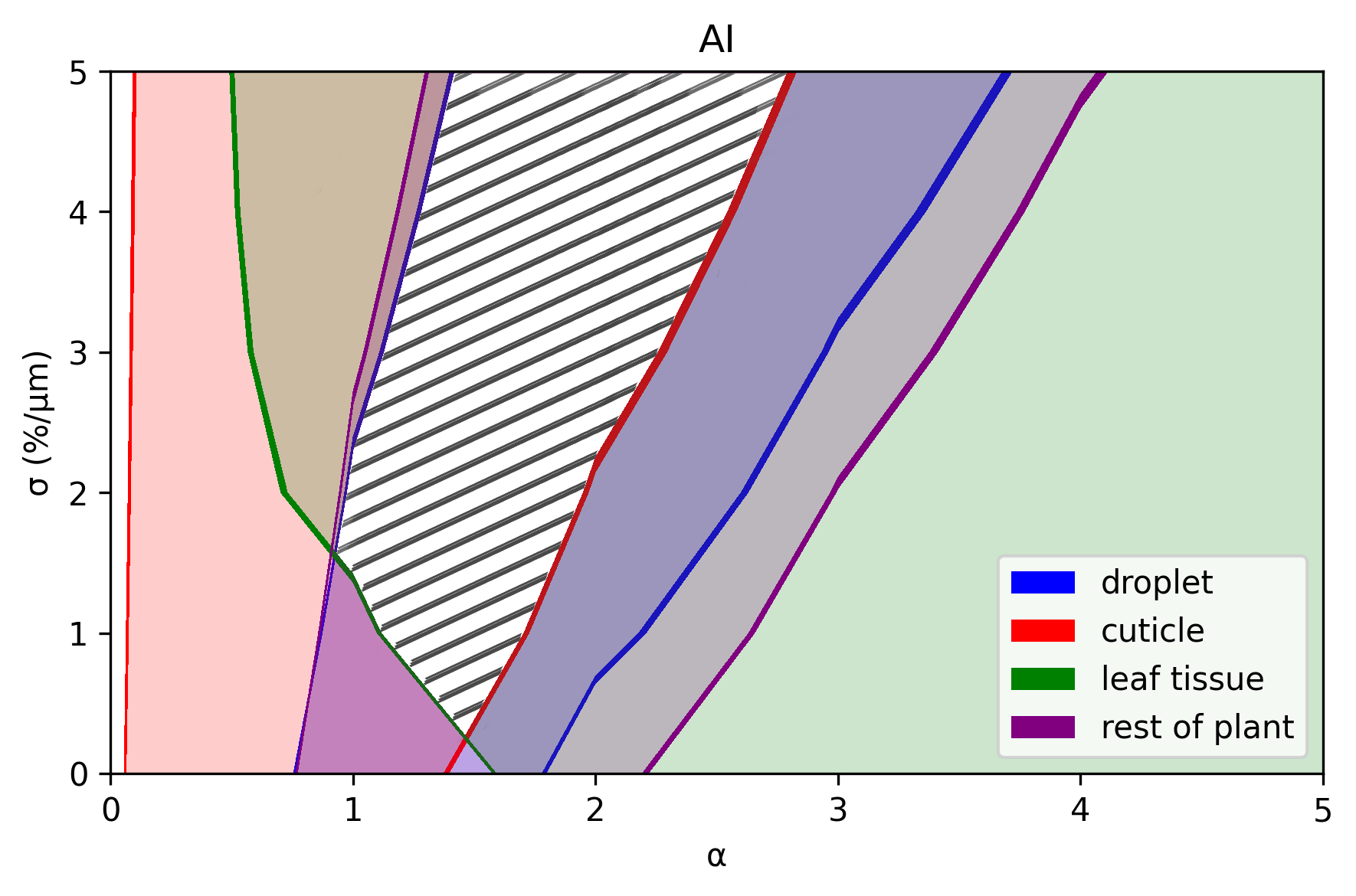}
			      \caption{Contour plots from \cref{fig: Numerical_Solution_alpha_sigma_DQ00.4} overlaid on one another to show ranges of $\alpha$ and $\sigma$ that produce the AI predicted data within the 95\% CIs of the experimental total amounts. The overlapping region is highlighted with black lines.}
			    \label{fig: Numerical_Solution_alpha_sigma_DQ00.4_overlay}
			\end{figure}

\begin{figure}[H]
			      \centering
				  \includegraphics[width=5cm]{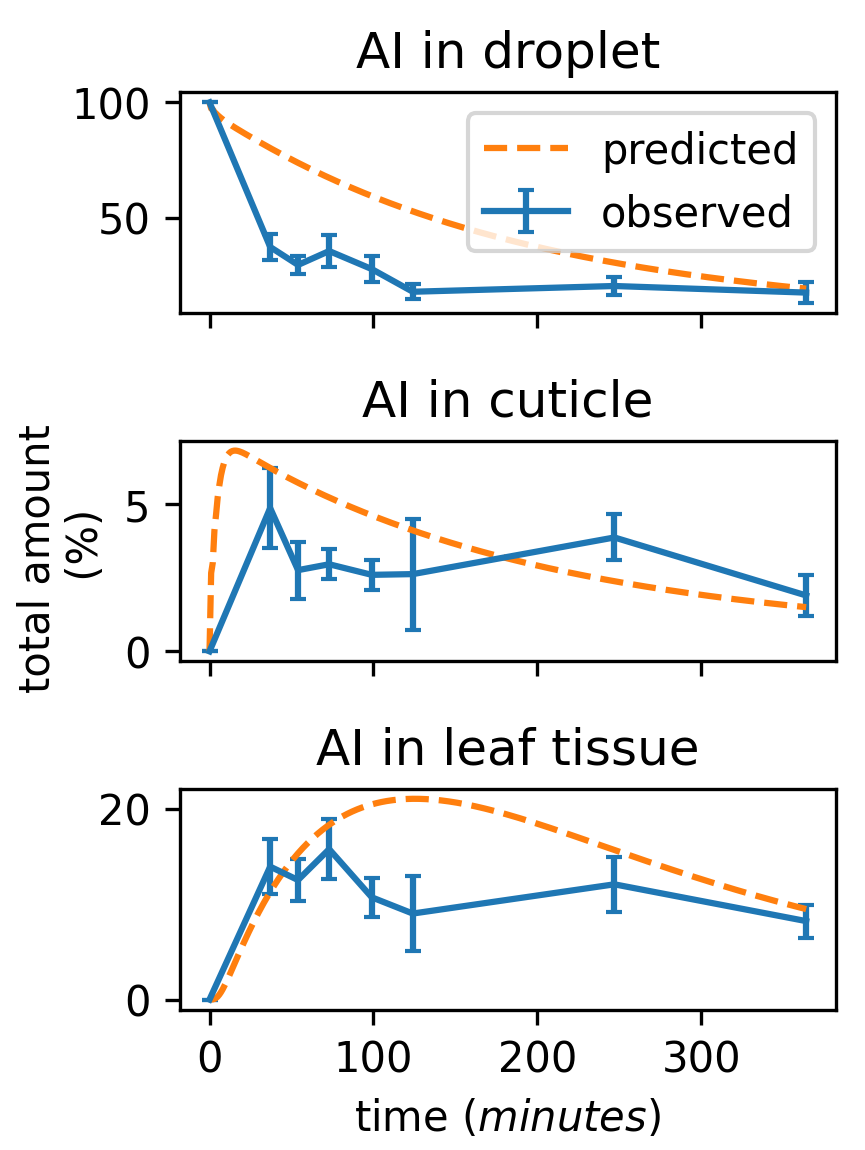}
				  \caption{Predicted uptake of AI using parameter estimates from \cref{subsec: Parameter Estimation} with $D_P,D_Q^0=0.4$ and $(\alpha,\sigma)=(1.5,3)$.}
			 	\label{fig: AI_subtotal_per_DP0.4_DQ00.4_alpha1.5_sigma3}
			\end{figure}

                We set $D_{Q,0}=0.4$, the lower bound from our parameter estimation and then allow $\alpha$ and $\sigma$ to vary. Results presented in \cref{fig: Numerical_Solution_alpha_sigma_DQ00.4} show how the change in these two parameters affects the amount of AI remaining in the droplet, cuticle, leaf tissue and rest of plant at the final time recorded in the experimental data. For each compartment, we superimpose the mean and 95\% CI from the experimental data which allows us to explore whether any pairs of values ($\alpha$, $\sigma$) might allow our model to give reasonable final time predictions. In \cref{fig: Numerical_Solution_alpha_sigma_DQ00.4_overlay} we overlay the plots from \cref{fig: Numerical_Solution_alpha_sigma_DQ00.4} to clearly highlight that a region in parameter space exists for which pairs of the parameters will give a reasonable fit to the end time data points. These values show that there is AJ influence on the AI diffusion since $\alpha>0$, and that saturation effect ($\sigma>0$) is present for this set of data. As an example, we take $(\alpha,\sigma)=(1.5,3)$ which is roughly central to the overlapping region in \cref{fig: Numerical_Solution_alpha_sigma_DQ00.4_overlay} and the corresponding fit is shown in \cref{fig: AI_subtotal_per_DP0.4_DQ00.4_alpha1.5_sigma3}.

    \section{Discussion}
	\label{sec: Discussion}
        We embarked on this work with two main goals: to create a mechanistic mathematical model to describe the uptake of pesticide through leaves and to include within that model the impact of the formulation on the transport of the AI. This was motivated by the growing momentum within the agro-chemical industry to exploit \textit{in-silico} modelling as a possible predictive tool as a way to reduce the costs, financial and environmental, on developing new products \cite{KalyabinaValeriyaP2021Pfdp}. 

        There is a relative paucity of existing models for leaf uptake (exceptions are, for example, \cite{Fujisawa2002,Rasmussen2003}), possibly in part due to the vast range of leaf structures, agro-chemical formulations and ambient environmental conditions \cite{WangC.J2007Fuop}. We see our work in this context as a novel contribution in particular through inclusion of the formulation effect on transport. 

        Working with a small data set, we used parsimonious assumptions to estimate key model parameters. Ideally, many of these parameters would be independently estimated using physico-chemical properties of the compounds and leaf environments. Our attempts to do this are presented in the \cref{sec: appendix_partition,sec: appendix_diffusion}. The estimates which we obtained using this approach did not result in good model fit, in particular because they produced high affinity for the compounds to remain in the leaf cuticle and movement across the cuticle was very slow. Further work in this area leading to a number of widely-applicable empirical relationships could ultimately reap significant reward because it would allow mechanistic modelling to replace some of the experimental work currently undertaken when developing new products. 

        The importance of including concentration-dependent diffusion for the AI is highlighted in \cref{fig: Numerical_Solution_alpha_sigma_DQ00.4_overlay} where we are able to identify a region in $(\alpha, \sigma)$ space where data from all compartments at the end time can effectively be represented. Since that region does not include $\alpha=0$, we conclude that the inclusion of concentration-dependent diffusion is required to observe this outcome. We are encouraged by the fit shown in \cref{fig: AI_subtotal_per_DP0.4_DQ00.4_alpha1.5_sigma3} that we are along the right direction in our modelling approach.

        Whether $\alpha=0$ or $\alpha>0$ (data not shown), the model consistently overestimates the amounts of both compounds in the leaf tissue at intermediate times. We anticipate that some of this effect would be moderated through the inclusion of other biological processes, metabolism or binding for example, known to occur in the leaf tissue. We chose not to include any such process in this work because the data is sufficiently coarse that we felt additional parameters to estimate would almost certainly lead to problems with parameter identifiability \cite{BEARUP2013,janzen2016}. Moreover, we were unable to find published data for either of these processes in the literature which would have helped with the estimation.

        There is plenty still to do before this model achieves its potential as a key tool in product development. However, we are heartened by our preliminary results which include identifying the need to look further into the physico-chemical empirical relationships for diffusion and partitioning. In the meantime, we will continue to work collaboratively with product development to build model capability.

\appendix

    \section{Partition coefficient estimates from literature}
	\label{sec: appendix_partition}
        Two estimates for the partition coefficient are presented in \cite{Rasmussen2003} which are both functions of the octanol-water partition coefficient, $P_{ow}$.
	\begin{subequations}
		\begin{align}
			\label{eq: par_coeff_1A}
			\log(K_{wax/water}) &= \log(P_{ow}) - 1, \\
			\label{eq: par_coeff_1B}
			\log(K_{cuticle/water}) &= -0.77+0.98 \log(P_{ow}),
		\end{align}	
	\end{subequations} 
        We use these estimates to obtain partition coefficients between wax and droplet, and between cuticle and leaf tissue given our assumption that the droplet and the leaf tissue can be described as aqueous environments. The results are given in \cref{tab: partition_coefficients_discussion}.
			
	\begin{table}[H]
		\centering
		\caption{Approximations for the partition coefficients using \cref{eq: par_coeff_1A} and \cref{eq: par_coeff_1B} from \cite{Rasmussen2003}.}
		\begin{tabular}{|c|c|c|c|}
			\hline
			&  Adjuvant &  Active Ingredient & unit \\ \hline
			$\log(P_{ow})$ & 3.90 & 3.19 & $-$ \\ \hline
			$\kappa_{1,A}$ & 797.33 & $-$ & $-$ \\ \hline
			$\kappa_{1,B}$ & 1127.20 & $-$ & $-$ \\ \hline
			$K_{1,A}$ & $-$ & 154.88 & $-$ \\ \hline
			$K_{1,B}$ & $-$ & 227.09 & $-$ \\ \hline
		\end{tabular}
		\label{tab: partition_coefficients_discussion}
	\end{table}
			
	\begin{figure}[H]
		\centering
		\includegraphics[width=11cm]{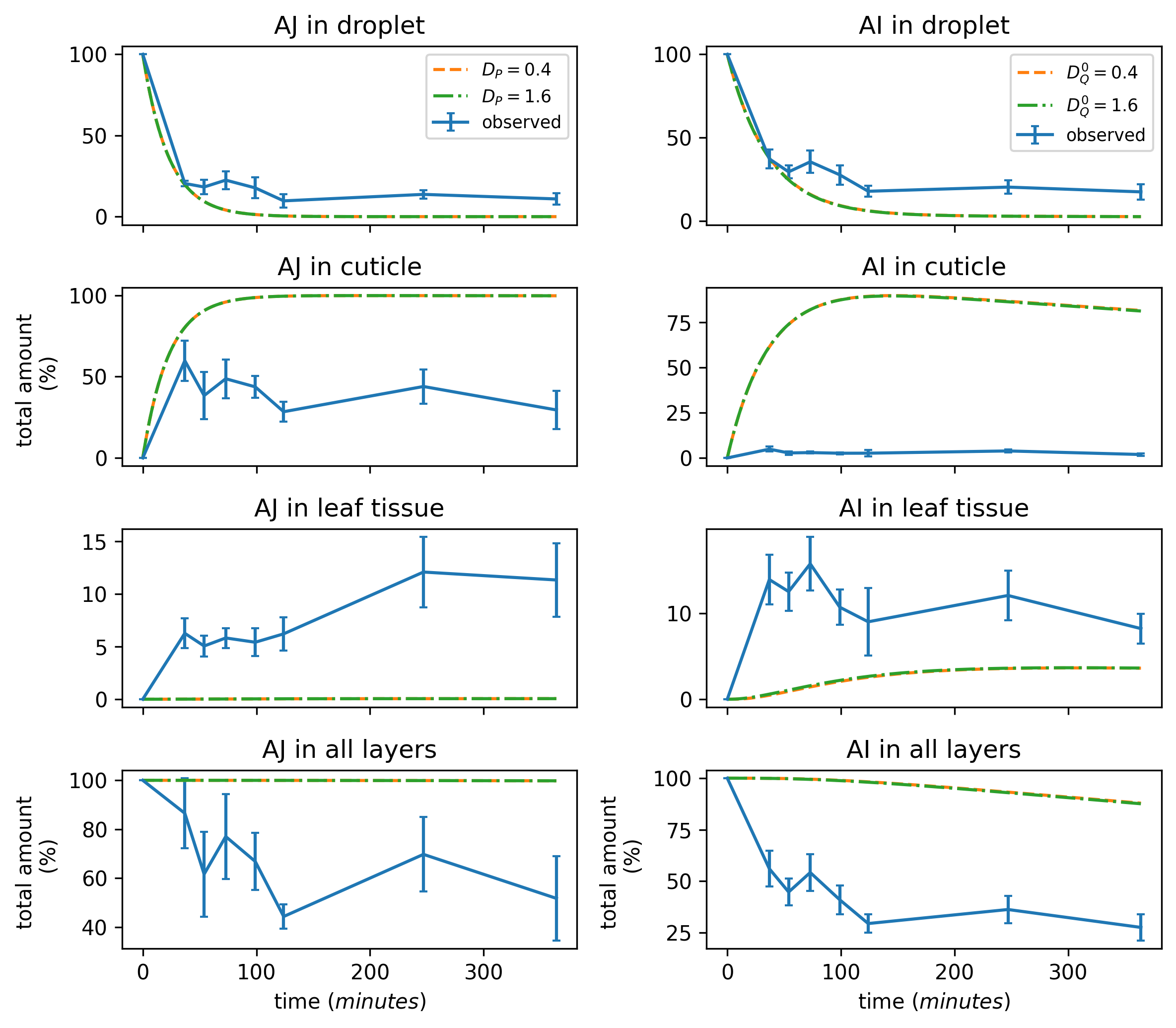}
		\caption{Predicted uptake of AJ and AI using parameter estimates from \cref{subsec: Parameter Estimation} with $\alpha=0$, and partition coefficients from \cref{tab: partition_coefficients_discussion}.}
		\label{fig: K_discussion}
	\end{figure}
        \Cref{fig: K_discussion} shows that using these estimates, transfer of both AI and AJ from cuticle to leaf tissue does not occur sufficiently and the model predicts that the chemical compounds get stuck in the cuticle. 

    \section{Diffusion coefficient estimates from literature}
	\label{sec: appendix_diffusion}
        Equations to estimate diffusion coefficients for both AJ and AI are also presented in \cite{Rasmussen2003} which are dependent on the McGowan volume \cite{Zhao2003} of a chemical compound. 
	\begin{subequations}
		\begin{align}
			\label{eq: diff_coeff_AJ}
			\log(D_{AJ}) &= -12.49-0.015 \times MV, \\
			\label{eq: diff_coeff_AI_wax}
			\log(D_{AI,wax}) &= -15.26 - 0.01 \times MV, \\
			\label{eq: diff_coeff_AI_cuticle}
			\log(D_{AI,cuticle}) &= -13 - 0.01 \times MV.
		\end{align}
	\end{subequations}
        We can see that there are two equations to estimate diffusion coefficient for the AI. This is due to the more complex structure of epicuticular wax on the outermost part of the wax/cuticle layer.
			
	\begin{table}[H]
		\centering
		\caption{Approximations for the diffusion coefficients using \cref{eq: diff_coeff_AJ,eq: diff_coeff_AI_wax,eq: diff_coeff_AI_cuticle} from \cite{Rasmussen2003}.}
		\begin{tabular}{|c|c|c|c|}
			\hline
			&  Adjuvant &  Active Ingredient & unit \\ \hline
			McGowan volume & 272.42 &  319.99 & ${cm}^{3}/mol$ \\ \hline
			$D_{AJ}$ & 2.65$\times10^{-17}$ & $-$ & ${m}^{2}/sec$ \\ \hline
			$D_{AI,wax}$ & $-$ &  3.47$\times10^{-19}$ & ${m}^{2}/sec$ \\ \hline
			$D_{AI,cuticle}$ & $-$ & 6.31$\times10^{-17}$ & ${m}^{2}/sec$ \\ \hline
			$D_{P}$ & 1.59$\times10^{-3}$ &  $-$ & ${\mu m}^2/min$ \\ \hline
			$D_{Q,0,wax}$ & $-$ &  2.08$\times10^{-5}$ & ${\mu m}^2/min$ \\ \hline
			$D_{Q,0,cuticle}$ & $-$ & 3.79$\times10^{-3}$ & ${\mu m}^2/min$ \\ \hline
		\end{tabular}
		\label{tab: diffusion_coefficients_discussion}
	\end{table}
        In this case, estimates of the diffusion coefficient mean that very little of either the AI or the AJ manage to cross the cuticle and again, the resulting predictions do not reflect the experimental data.

    \section*{Acknowledgments}
        We would like to thank Mr. Damon Raziel Salvatore for his help with producing the diagrams used in this paper.

    \bibliographystyle{siamplain}
    \bibliography{References}

\end{document}